\documentclass[11pt,a4paper,english]{article}

\usepackage{verbatim}
\usepackage{amsfonts,amsmath,amsthm}

\def\LaTeX{L\kern -.36em\raise .3ex\hbox{\sc a}\kern -.15em T\kern -.1667em%
\lower .7ex\hbox\mathbb{E}\kern -.125em X}

\usepackage{amssymb}
\usepackage{float}
\usepackage{epsfig}
\usepackage{amsmath}
\usepackage[english]{babel}
\usepackage{amstext}
\usepackage{mathrsfs}
\usepackage{amsthm}
\usepackage{fancyhdr}
\usepackage{latexsym}
\usepackage{amsmath}
\usepackage{amsfonts}
\usepackage{dsfont}
\usepackage{amssymb}
\usepackage{fancyhdr}

\usepackage[usenames]{color}
\definecolor{red}{rgb}{1.0,0.0,0.0}

\definecolor{blu}{rgb}{0.0,0.0,1.0}

\def\F{{\cal F}}

\def\to{\longrightarrow}

\def \0{{\textbf{0}}}
\def\norm{{\| \kern -.05em | }}

\newtheorem{Theorem}{Theorem}[section]
\newtheorem{Lemma}[Theorem]{Lemma}
\newtheorem{Corollary}[Theorem]{Corollary}
\newtheorem{Proposition}[Theorem]{Proposition}
\newtheorem{Definition}[Theorem]{Definition}
\newtheorem{Remark}[Theorem]{Remark}
\newtheorem{Assumption}[Theorem]{Assumption}

\renewcommand{\theequation}{\thesection.\arabic{equation}}

\def \N{\mathbb{N}}
\def \R{\mathbb{R}}

\def \E{\mathbb{E}}
\def \F{\mathbb{F}}

\def \P{\mathbb{P}}

\def \Ac{{\cal A}}

\def \Cc{{\cal C}}
\def \Dc{{\cal D}}

\def \Fc{{\cal F}}

\def \Ic{{\cal I}}
\def \Kc{{\cal K}}

\def \Oc{{\cal O}}

\def \Tc{{\cal T}}

\def \eps{\varepsilon}

\def \ep{\hbox{ }\hfill$\Box$}

\def\beqs{\begin{eqnarray*}}
\def\enqs{\end{eqnarray*}}
\def\beq{\begin{eqnarray}}
\def\enq{\end{eqnarray}}

\begin{document}

\markright{ \sc Optimal Boundary of an Irreversible Investment Problem}

\makeatletter  
\@addtoreset{equation}{section}
\def\theequation{\arabic{section}.\arabic{equation}}
\makeatother  

\title{\textbf{Optimal Boundary Surface for Irreversible Investment with Stochastic Costs}\footnote{The first author was supported by EPSRC grant EP/K00557X/1; Financial support by the German Research Foundation (DFG) via grant Ri--1128--4--2 is gratefully acknowledged by the third author. This work was started during a visit of the second author at the Center for Mathematical Economics (IMW) at Bielefeld University thanks to a grant by the German Academic Exchange Service (DAAD). The second author thankfully aknowledges the financial support by DAAD and the hospitality of IMW.} }
\author{Tiziano De Angelis\thanks{School of Mathematics, University of Leeds, Woodhouse Lane, Leeds LS2 9JT, United Kingdom; \texttt{t.deangelis@leeds.ac.uk}}
\and Salvatore Federico\thanks{Dipartimento di Scienze per l'Economia e l'Impresa, Universit\`{a} degli Studi di Firenze, Via delle Pandette 9, 50127 Firenze, Italy; \texttt{salvatore.federico@unifi.it}}
\and Giorgio Ferrari\thanks{Center for Mathematical Economics (IMW), Bielefeld University, Universit\"atsstrasse 25, D-33615 Bielefeld, Germany; \texttt{giorgio.ferrari@uni-bielefeld.de}}}
\date{\today}
\maketitle

\vspace{0.5cm}

{\textbf{Abstract.}} This paper examines a Markovian model for the optimal irreversible investment problem of a firm aiming at minimizing total expected costs of production. We model market uncertainty and the cost of investment per unit of production capacity as two independent one-dimensional regular diffusions, and we consider a general convex running cost function. The optimization problem is set as a three-dimensional degenerate singular stochastic control problem. We provide the optimal control as the solution of a reflected diffusion at a suitable boundary surface. Such boundary arises from the analysis of a family of two-dimensional parameter-dependent optimal stopping problems and it is characterized in terms of the family of unique continuous solutions to parameter-dependent nonlinear integral equations of Fredholm type.
\smallskip


\smallskip


\smallskip



\section{Introduction}

In this paper we study a Markovian model for a firm's optimal irreversible investment problem. The firm aims at minimizing total expected costs of production when its running cost function depends on the uncertain condition of the economy as well as on on the installed production capacity, and the cost of investment per unit of production capacity is random. In mathematical terms, this amounts to solving the three-dimensional degenerate singular stochastic control problem
\beq
\label{OCPintro}
V(x,y,z):=\inf_{\nu}\mathbb{E}\bigg[\int_0^{\infty} e^{-r t} c(X^x_t, z + \nu_t){dt}+ \int_0^{\infty} e^{-r t} Y^y_t d \nu_t \bigg],
\enq
where the infimum is taken over a suitable set of nondecreasing admissible controls. Here $X$ and $Y$ are two 
diffusion processes modeling market uncertainty and the cost of investment per unit of production capacity, respectively. The control process $\nu_t$ is the cumulative investment made up to time $t$ and $c$ is a general convex cost function. We solve problem \eqref{OCPintro} by relying on the connection existing between singular stochastic control and optimal stopping (see, e.g., \cite{BK}, \cite{BatherChernoff}, \cite{Benth}, \cite{BoetiusKohlmann}, \cite{KS1} and \cite{KW}). In fact, we provide the optimal investment strategy $\nu^*$ in terms of an optimal boundary surface $(x,y)\mapsto z^*(x,y)$ that splits the state space into \emph{action} and \emph{inaction} regions. Such surface 
is then uniquely characterized through a family of continuous solutions to parameter-dependent, nonlinear integral equations of Fredholm type.

In the mathematical economic literature, singular stochastic control problems are often employed to model the irreversible (partially reversible) optimal investment problem of a firm operating in an uncertain environment (see \cite{CH1}, \cite{DeAF2013}, \cite{FedericoPham}, \cite{freeandbase}, \cite{GP}, \cite{Kobila}, \cite{MZ}, \cite{RS} and references therein, among many others). The monotone (bounded-variation) control represents the cumulative investment (investment-disinvestment) policy used by the firm to maximize total net expected profits or, alternatively, minimize total expected costs. The optimal timing problem, associated to the optimal investment one, is then related to real options, as pointed out by \cite{McDonaldSiegel} and \cite{Pindyck} among others.

Problems of stochastic irreversible (or partially reversible) investment can be tackled via a number of different approaches. Among others, these include dynamic programming techniques (see, e.g., \cite{FedericoPham}, \cite{GP}, \cite{Kobila} and \cite{MZ}), stochastic first-order conditions and the Bank-El Karoui's Representation Theorem \cite{BankElKaroui} (see, e.g., \cite{B}, \cite{CF}, \cite{freeandbase} and \cite{RS}), the transformation method of \cite{BayEgami} in the case of one-dimensional problems, and the analytical study of non linear PDEs with gradient constraints (see for example \cite{SS90} and \cite{SS91}). The introduction of a stochastic investment cost $Y_t$ is very natural from the point of view of economic modelling (see e.g.~\cite{B}); nevertheless, it makes the analysis of the optimal boundary of \eqref{OCPintro} rather difficult. 

The three-dimensional structure of our problem \eqref{OCPintro} makes seemingly hopeless a direct study of the associated Hamilton-Jacobi-Bellman equation with the aim of finding explicit smooth solutions (as in the two-dimensional problem of \cite{MZ}, among others). In fact, in our case the linear part of the Hamilton-Jacobi-Bellman equation for the value function of problem \eqref{OCPintro} is a PDE (rather than a ODE) and it does not have a general solution.
On the other hand, arguing as in \cite{freeandbase}, we might tackle problem \eqref{OCPintro} by relying on a stochastic first-order conditions approach and a suitable application of the Bank-El Karoui's representation theorem \cite{BankElKaroui}.
However, the integral equation for the optimal boundary, which derives from the main result of \cite{freeandbase} (i.e., \cite[Th.~3.11]{freeandbase}), cannot be obtained in our multi-dimensional setting.

In this paper we study problem \eqref{OCPintro} by relying on the connection between singular stochastic control and optimal stopping. Building on a well known result concerning variational inequalities (see Proposition \ref{teo:ex}), we then develop almost exclusively probabilistic arguments to find an optimal control $\nu^*$. We show that such $\nu^*$ is the minimal effort needed to keep the (optimally controlled) state process above an optimal boundary surface $z^*$, whose level curves $z^*(x,y)=z$, where $z\in\R^+$, are the optimal boundaries $x\mapsto y^*(x;z)$ of the parameter-dependent optimal stopping problems associated to the original singular control one. Under some further mild conditions, we characterize each function $y^*(\,\cdot\,;z)$, $z \in\R^+$, as the unique continuous solution of a nonlinear integral equation of Fredholm type (see our Theorem \ref{teo:uniqueness} below). 

One should notice that the connection to optimal stopping was also used, for example, in \cite{SS90} to study a non linear PDE problem with gradient constraint related to an optimization similar to \eqref{OCPintro}. In \cite{SS90} $X\in\R^{n-1}$ is a Brownian motion, the controlled process $Z\in \R$ is a linearly controlled Brownian motion, $Y_\cdot\equiv1$, and a detailed analysis of the optimal boundary (the free-boundary of the PDE) is carried out based exclusively on analytical methods. An extension of those methods to our setting seems possible, but non trivial. Here, instead we develop a different approach mostly employing stochastic calculus to uniquely determine our optimal boundary.

The issue of finding integral equations for the optimal boundary of optimal stopping problems has been successfully addressed in a number of papers (cf.~\cite{PeskShir} for a survey) and dates back to the work of Van Moerbeke \cite{VM} among others (see also \cite{Cannon} for a survey of PDE methods). In the context of one-dimensional stochastic (ir)reversible investment problems on a finite time-horizon integral equations for the optimal boundaries have been obtained recently by an application of Peskir's local time-space calculus (see \cite{DeAF2013} and references therein for details). However, those arguments cannot be applied in our case, since it seems quite hard to prove that the process $\{y^*(X^x_t;z),\,t\ge0\}$ is a semimartingale for each given $z \in \R^+$, as it is required in \cite[Th.\ 2.1]{Peskir}.
On the other hand, numerically computable integral equations for multi-dimensional settings have been studied, for instance in \cite[Sec.~13]{PeskShir}, where a diffusion $X$ was considered along with its running supremum $S$. However, unlike \cite[Sec.~13]{PeskShir} here we deal with a genuine two dimensional diffusion $(X,Y)$ with $X$ and $Y$ independent. This gives rise to a completely different analysis of the problem and new methods have been developed.

In summary, the main contributions of our work are the following: $i)$ we provide optimal boundaries for models of irreversible investment under uncertainty with stochastic investment costs; $ii)$ as a byproduct we develop methods to uniquely characterize optimal boundaries of infinite time horizon optimal stopping problems for 2-dimensional diffusions, thus extending portions of the existing techniques based on stochastic calculus.
These optimal boundaries might also be numerically treated relying on numerical methods for nonlinear Fredholm integral equations of second kind (see Remark \ref{rem-inteq} below).

The paper is organized as follows. In Section \ref{sec:IrrInvProb} we set the stochastic irreversible investment problem. In Sections \ref{sec:associated} and  \ref{sec:char} we introduce the associated family of optimal stopping problems and we characterize its value functions and its optimal-boundaries. The form of the optimal control is provided in Section \ref{sec:opt}. Finally, some technical results are discussed in Appendix \ref{app:os1d}.


\section{The Stochastic Irreversible Investment Problem}
\label{sec:IrrInvProb}

In this section we set the stochastic irreversible investment problem object of our study. Let $(\Omega, \mathcal{F},(\Fc_t)_{t\ge0}, \mathbb{P})$ be a complete filtered probability space with $\F=\{\mathcal{F}_t, t \geq 0\}$ the filtration generated by a two-dimensional Brownian motion $W=\{�(W^1_t,W^2_t),\, t\geq 0\}$ and augmented with $\P$-null sets.

\begin{enumerate}
\item A real process $X=\{X_t,\, t\geq 0\}$ represents the uncertain status of the economy (typically, the demand of a good or, more generally, some indicator of macroeconomic conditions). We assume that $X$ is a time-homogeneous Markov diffusion satisfying the stochastic differential equation (SDE)
\beq\label{state:X}
dX_t=\mu_1(X_t)dt+\sigma_1(X_t)dW_t^1,\qquad X_0=x,
\enq
for some Borel functions $\mu_1$ and $\sigma_1$ to be specified. To account for the dependence of $X$ on its initial position we denote the solution of \eqref{state:X} by $X^x$.
\item A one-dimensional positive process $Y=\{Y_t,\, t\geq 0\}$ represents the cost of investment per unit of production capacity. We assume that $Y$ evolves according to the SDE
\beq\label{state:Y}
dY_t=\mu_2(Y_t)dt+\sigma_2(Y_t)dW_t^2,\qquad Y_0=y,
\enq
for some Borel functions $\mu_2$ and $\sigma_2$ to be specified as well. Again, to account for the dependence of $Y$ on $y$, we denote the solution of \eqref{state:Y} by $Y^y$.
\item A control process $\nu=\{\nu_t,\ t\geq 0\}$ describes an investment policy of the firm and $\nu_t$ is the cumulative investment made up to time $t$. We say that a control process $\nu$ is admissible if it belongs to the nonempty convex set
\begin{align}\label{setadmissblecontrols}
\mathcal{V}:=\{ \nu:\Omega\times \mathbb{R}^+ \mapsto \mathbb{R}^+\,|\, \mbox{$t\mapsto\nu_t$ is c\`adl\`ag, nondecreasing, $\F$-adapted}\}.
\end{align}
In the following we set $\nu_{0^-} =0$, for every $\nu\in\mathcal{V}.$
\item A purely controlled process $Z=\{Z_t,\,t\geq 0\}$, represents the production capacity of the firm and it is defined by
\beq
\label{state:Z}
Z_t:=z+\nu_t, \qquad z \in \R^+.
\enq
The process $Z$ depends on its initial position $z$ and on the control (investment) process $\nu$, therefore we denote it by $Z^{z,\nu}$.

\end{enumerate}
We assume that the uncontrolled diffusions $X^x$ and $Y^y$ have state-space $\mathcal{I}_1=(\underline{x},\overline{x}) \subseteq \R$ and $\mathcal{I}_2=(\underline{y},\overline{y}) \subseteq \R^+$, respectively. 
\begin{Remark}
\label{rem:ass ind}
Some of our results would continue to hold also if $X\in\mathbb{R}^n$ and the Brownian motions driving the SDEs for $X$ and $Y$ were correlated. Indeed, the only proofs employing independence of $X$ and $Y$, and $X\in\R$ are those of Proposition \ref{prop:y*}, Proposition \ref{continuousboundary} and Theorem \ref{teo:uniqueness}. However, since these are key results in our paper, we adopt the above setting from the beginning to simplify the exposition.
\end{Remark}
The boundary behaviour of $X^x$ and $Y^y$ and further requirements on the coefficients $\mu_i,\sigma_i$, $i=1,2$ are specified in the following assumption.

\begin{Assumption}\label{ass:D2}\hspace{10cm}
\begin{itemize}

\item[(i)]
The coefficients $\mu_i: \R \mapsto \R,\ \sigma_i:\R \mapsto \R^+ $, $i=1,2$, are such that
$$|\mu_i(\zeta)-\mu_i(\zeta')|\ \leq \ K_i|\zeta-\zeta'|, \ \ \ |\sigma_i(\zeta)-\sigma_i(\zeta')|\ \leq \ M_i|\zeta-\zeta'|^{\gamma}, \ \ \ \forall \zeta,\zeta'\in \mathcal{I}_i,$$
for some $K_i>0$, $M_i>0$ and $\gamma \in [\frac{1}{2},1]$.
\item[(ii)]
The diffusions $X^x$ and $Y^y$ are nondegenerate, i.e.\ $\sigma^2_i>0$ in $\mathcal{I}_i$, $i=1,\,2$.
\item[(iii)]
The boundaries $\underline{x},\overline{x}$ are non-exit for the diffusion $X^x$ and the boundaries $\underline{y},\overline{y}$ are natural for the diffusion $Y^y$.\footnote{A boundary point $\xi$ is non-exit for a diffusion process if: \emph{i)} the process started from the interior of its state space cannot reach $\xi$ in finite time, and \emph{ii)} the process starting from $\xi$ immediately enters the interior of the state space. On the other hand, a non-entrance boundary point $\xi$ can be reached in a  finite time but it cannot be a starting point for the diffusion. Finally, a boundary point $\xi$ is natural if it is: non-entrance and non-exit (cf.~for instance \cite[Ch.~2, p.~15]{BS}). Moreover, for $\xi$ natural and finite one also has $\mu_2(\xi)=\sigma_2(\xi)=0$ if $\xi=\underline{y}$ (or $\xi=\overline{y}$). That is shown in Appendix \ref{A1} for the sake of completeness.}
\end{itemize}
\end{Assumption}

Assumption \ref{ass:D2} guarantees that
\beq
\label{LI}
\int_{\zeta-\varepsilon_o}^{\zeta+\varepsilon_o}\frac{1 + |\mu_i(y)|}{|\sigma_i(y)|^2}\,dy < +\infty,\quad\text{for some $\varepsilon_o>0$ and every $\zeta$ in $\mathcal{I}_i$},
\enq
 hence both \eqref{state:X} and \eqref{state:Y} have a weak solution that is unique in the sense of probability law (cf.\ \cite[Ch.\,5.5]{KS}). Such solutions do not explode in finite time, due to the fact that the coefficients have at most linear growth. On the other hand, Assumption \ref{ass:D2}-$(i)$ also guarantees pathwise uniqueness for the solutions of \eqref{state:X} and \eqref{state:Y} by the Yamada-Watanabe result (see \cite[Ch.\ 5.2, Prop.\ 2.13]{KS} and \cite[Ch.\ 5.3, Rem.\ 3.3]{KS}, among others). Therefore, \eqref{state:X} and \eqref{state:Y} have a unique strong solution due to \cite[Ch.\ 5.3, Cor.\ 3.23]{KS} for any $x\in\mathcal{I}_1$ and $y\in\mathcal{I}_2$.
Also, it follows from (\ref{LI}) that the diffusion processes $X^{x}$ and $Y^y$ are regular in $\mathcal{I}_1$ and $\mathcal{I}_2$, respectively; that is, $X^{x}$ (resp., $Y^y$) hits a point $\zeta$ (resp., $\zeta'$) with positive probability, for any $x$ and $\zeta$ in $\mathcal{I}_1$ (resp., $y$ and $\zeta'$ in $\mathcal{I}_2$). Hence the state spaces $\mathcal{I}_1$ and $\mathcal{I}_2$ cannot be decomposed into smaller sets from which $X^{x}$ and $Y^y$ could not exit (see, e.g., \cite[Ch.\ V.7]{Rog}).
Finally, there exist continuous versions of $X^{x}$ and $Y^y$ and we shall always refer to those versions throughout this paper.

Assumption \ref{ass:D2} implies the comparison criterion (see, e.g., \cite[Ch.\,5.2, Prop.\,2.18]{KS}); i.e.,\,
\beq\label{comparison}
x,\,x'\in\Ic_1\,,\quad x \leq x'  \Longrightarrow  X_{t}^{x} \leq X_{t}^{x'}, \quad\text{$\P$-a.s.~$\forall t\ge0$.}
\enq
Moreover, repeating arguments as in the proof of \cite[Ch.\,5.2, Prop.\,2.13]{KS} one also finds
\beq\label{contdata}
x_n \rightarrow\, x_0 \  \mbox{in $\Ic_1$ as} \  n\rightarrow\infty \ \Longrightarrow \ X^{x_n}_t\stackrel{L^1}{\to} X_t^{x_0} \ \Longrightarrow \ X^{x_n}_t\stackrel{\P}{\to} X_t^{x_0}, \ \  \ \forall t\geq 0;
\enq
Analogously, for the unique solution of \eqref{state:Y} one has
\beq\label{comparisonY}
y,\,y'\in\Ic_2\,,\quad y \leq y'  \Longrightarrow  Y_{t}^{y} \leq Y_{t}^{y'}, \quad\text{$\P$-a.s.~$\forall t\ge0$};
\enq
and
\beq\label{contdataY}
y_n \rightarrow\, y_0 \  \mbox{in $\Ic_2$ as} \  n\rightarrow\infty \ \Longrightarrow \ Y^{y_n}_t\stackrel{L^1}{\to} Y_t^{y_0}\Longrightarrow \ Y^{y_n}_t\stackrel{\P}{\to} Y_t^{y_0}, \ \  \ \forall t\geq 0.
\enq
Standard estimates on the solution of SDEs with coefficients having sublinear growth imply that (cf., e.g., \cite[Ch.\ 2.5, Cor.\ 12]{K})
\beq
\label{stimaKrylovXY}
\mathbb{E}\Big[|X^x_t|^q\Big] \leq \kappa_{0,q}(1 + |x|^q)e^{\kappa_{1,q}t}, \qquad
\mathbb{E}\Big[|Y^y_t|^q\Big] \leq \theta_{0,q}(1 + |y|^q)e^{\theta_{1,q}t}, \qquad t\geq 0,
\enq
for any $q \geq 0$, and for some $\kappa_{i,q}:=\kappa_{i,q}(\mu_1,\sigma_1)>0$ and $\theta_{i,q}:=\theta_{i,q}(\mu_2,\sigma_2)>0$, $i=0,1$.

Within this setting we consider a firm that incurs investment costs and a running cost $c(x,z)$ depending on the state of economy $x$ and the production capacity $z$. The firm's total expected cost of production associated to an investment strategy $\nu\in\mathcal{V}$ is
\beq \label{objective}
\mathcal{J}_{x,y,z}(\nu):=\mathbb{E}\bigg[\int_0^{\infty} e^{-r t} c(X^x_t,Z^{z,\nu}_t){dt}+ \int_0^{\infty} e^{-r t} Y^y_t d \nu_t \bigg],
\enq
for any $(x,y,z) \in \mathcal{I}_1 \times \mathcal{I}_2 \times \R^+$. Here $r$ is a positive discount factor and the cost function satisfies
\begin{Assumption}
\label{ass:cost} \hspace{10cm}
\begin{itemize}

\item[(i)]
$c: \overline{\mathcal{I}}_1 \times \mathbb{R}^+ \mapsto \mathbb{R}^+$ is such that $c\in C^0(\overline{\mathcal{I}}_1\times \R^+;\R^+)$, $c(x,\cdot)\in C^1(\mathbb{R}^+)$ for every $x\in\overline{\mathcal{I}}_1$, and  $c_z\in C^{\alpha}(\overline{\mathcal{I}}_1\times \R^+; \R)$ for some $\alpha > 0$ (that is, $c_z$ is $\alpha$-H\"older continuous).

\item[(ii)]  $c(x,\cdot)$ is convex for all $x\in\overline{\mathcal{I}}_1$ and $c_z(\cdot,z)$ is nonincreasing for every $z\in\mathbb{R}^+$.

\item[(iii)] $c$ and $c_z$ satisfy a polynomial growth condition with respect to x; that is, there exist locally bounded functions $\eta_o, \gamma_o: \R^+ \mapsto \R^+$, and a constant $\beta \geq 0$ such that
$$|c(x,z)| + |c_z(x,z)| \leq \eta_o(z) + \gamma_o(z)|x|^{\beta}.$$
\end{itemize}
\end{Assumption}

Throughout this paper we also make the following standard assumption that guarantees in particular finiteness for our problem (see Remark \ref{rem:c}-(\textbf{3}) and Lemma \ref{propertiesV} below)
\begin{Assumption}
\label{Ass:r}
$r > \kappa_{1,\beta},$
\end{Assumption}
\noindent with $\kappa_{1,q}$ as in \eqref{stimaKrylovXY} and with $\beta$ of Assumption \ref{ass:cost}-$(iii)$.

\begin{Remark}\label{rem:c}
\noindent {\bf 1.} Any function $c$ of the spread $|x-z|$ between capacity and demand in the form
\beq\label{gspread}
c(x,z)= K_0 |x-z|^{\delta}, \ \ \ K_0\geq 0, \ \delta>1,
\enq
satisfies Assumption \ref{ass:cost}. We observe that \eqref{gspread} is a natural choice, e.g., in an energy market framework where $x$ represents the demand net of renewables (thus having stochastic nature) and $z$ the amount of conventional supply. Failing to meet the demand as well as an excess of supply generate costs for the energy provider.

\noindent {\bf 2.} The second part of Assumption \ref{ass:cost}-(ii) captures the negative impact on marginal costs due to an increase of demand. It is intuitive in \eqref{gspread} that an increase of $z$ will produce a reduction (increase) of costs which is more significant the more the demand is above (below) the supply.

\noindent {\bf 3.} It follows from \eqref{stimaKrylovXY}, Assumption \ref{ass:cost}-(iii) and Assumption \ref{Ass:r} that $c$ and $c_z$ satisfy the integrability conditions
\begin{itemize}
\item[(a)]
$\displaystyle \mathbb{E}\bigg[\int_0^\infty e^{-rt} c(X_t^x,z)dt\bigg]<\infty, \ \ \forall (x,z) \in \mathcal{I}_1\times\R^+$;
\item[(b)]
$\displaystyle \mathbb{E}\bigg[\int_0^\infty e^{-rt} |c_z(X_t^x,z)|dt\bigg]<\infty, \ \ \forall (x,z)\in \mathcal{I}_1\times\mathbb{R}^+.$
\end{itemize}
Notice that integrability property (b) above guarantees finiteness of the value of the optimal stopping problem we will discuss in the next section (see \eqref{P:OS}).

\noindent {\bf 4.} In the benchmark case of $X^x$ given by a geometric Brownian motion with drift $\mu_1$ and volatility $\sigma_1$, using the well known formula for the Laplace transform of a Gaussian random variable, one obtains
$$\mathbb{E}[(X^x_t)^{\beta}] = x^{\beta}\exp\Big\{ \beta \Big(\mu_1 + \frac{1}{2}\sigma_1^2(\beta-1)\Big)t\Big\}$$
and therefore Assumption 2.3 reads $r> \beta(\mu_1 + \frac{1}{2}\sigma_1^2(\beta-1))=:\kappa_{1,\beta}.$

\noindent {\bf 5.} It is worth noticing that all the results of this paper hold even if we allow the running cost function $c$ depending on the triple $(x,y,z)$, satisfying conditions analogous to (i)-(iii) of Assumption \ref{ass:cost} and with $y \mapsto c_z(x,y,z)$ increasing.
However, since this extension does not have a clear economic meaning and in order to simplify the exposition, we only consider $c$ as in Assumption \ref{ass:cost} above.

\end{Remark}

The firm's manager aims at picking an irreversible investment policy $\nu^* \in \mathcal{V}$ (cf.~\eqref{setadmissblecontrols}) that minimizes the total expected cost \eqref{objective}.
Therefore, by denoting the state space $\mathcal{O}:=\mathcal{I}_1\times \mathcal{I}_2\times\mathbb{R}^+,$ the firm's manager is faced with the optimal irreversible investment problem with value function
\begin{equation}\label{OCP}
V(x,y,z):=\inf_{\nu\in\mathcal{V}}\mathcal{J}_{x,y,z}(\nu), \ \ \ (x,y,z)\in\mathcal{O}.
\end{equation}

Notice that \eqref{stimaKrylovXY}, Assumption \ref{ass:cost} and Assumption \ref{Ass:r} (cf.\ also Remark \ref{rem:c}-(\textbf{3})), together with the affine nature of $Z^{z,\nu}$ in the control variable lead to the following
\begin{Lemma}
\label{propertiesV}
The value function $V(x,y,z)$ of \eqref{OCP} is finite for all $(x,y,z)\in\mathcal{O}$ and such that $z \mapsto V(x,y,z)$ is convex.
\end{Lemma}

Problem \eqref{OCP} is a degenerate, three-dimensional, convex singular stochastic control problem of monotone follower type (see, e.g., \cite{KaratzasElKarouiSkorokhod}, \cite{KS1} and references therein). Moreover, if $c(x,\cdot)$ is strictly convex, then $\mathcal{J}_{x,y,z}(\cdot)$ of (\ref{objective}) is strictly convex on $\mathcal{V}$ as well, and hence if a solution to (\ref{OCP}) exists, it must be unique. 



\section{The Family of Associated Optimal Stopping Problems}
\label{sec:associated}


We now introduce and study the family of optimal stopping problems that we expect to be associated to the singular control problem \eqref{OCP} (see \cite{BK} among others).
Set $$\mathcal{T}:=\{\tau\,:\,\tau\:\text{are $\mathbb{F}$-stopping times}\},$$
and define
\begin{align}\label{functional:OS}
\Psi_{x,y,z}(\tau):= \mathbb{E}\bigg[\int_0^\tau e^{-rt} c_z(X_t^x,z)dt-e^{-r\tau}Y^y_\tau \bigg], \quad \tau \in \mathcal{T},  \quad (x,y)\in \mathcal{I}_1\times \mathcal{I}_2, \ \ z\in\mathbb{R}^+.
\end{align}
For any $z\in\mathbb{R}^+$ we consider the optimal stopping problem
\begin{align}\label{P:OS}
v(x,y;z):= \sup_{\tau\in \mathcal{T}} \Psi_{x,y,z}(\tau), \qquad (x,y)\in \mathcal{I}_1\times \mathcal{I}_2
\end{align}
and notice that $\big\{v(x,y;z),\,z\in\mathbb{R}^+\big\}$ is a family of two-dimensional parameter-dependent optimal stopping problems.

In the rest of the present section and in the next one, we fix $z \in \R^+$ and we study the optimal stopping problem \eqref{P:OS}. Denote its state space by $Q:=\mathcal{I}_1\times\mathcal{I}_2$. We introduce the following (cf.\ \cite[Ch.\ 1, Def.\ 4.8]{KS})
\begin{Definition}
\label{classD}
A right-continuous stochastic process $\xi:=\{\xi_t, t \geq 0\}$ is of class (D) if the family of random variables $\{\xi_{\tau}\mathds{1}_{\{\tau < \infty\}},  \tau \in \mathcal{T}\}$ is uniformly integrable,
\end{Definition}
\noindent and we make the following technical
\begin{Assumption}\label{ass:ui}
The process $\{e^{-rt}Y^y_t, t \geq 0\}$ is of class (D) and such that $\lim_{t \rightarrow \infty}e^{-rt}Y^y_t = 0$ $\P$-a.s.
\end{Assumption}

\begin{Remark}
\label{Rem:ass:ui}
\noindent {\bf 1.} The  process $e^{-rt}Y^y_t$ is of class (D) if, e.g., $\mathbb{E}[\sup_{t\geq 0}e^{-rt}Y^y_t]< \infty$, a standard technical assumption in the general theory of optimal stopping (see, e.g, \cite[Ch.\ I]{PeskShir}).

\noindent {\bf 2.} The last requirement of Assumptions \ref{ass:ui} is satisfied if, e.g., $\{e^{-rt}Y^y_t, \,t \geq 0\}$ is an ($\mathcal{F}_t$)-supermartingale and $r>\theta_{1,1}$ (cf.\ \eqref{stimaKrylovXY}). Indeed from \cite[Ch.\ 1, Problem 3.16]{KS} and Fatou's Lemma one has
$$\displaystyle 0  \leq \mathbb{E}[\lim_{t \rightarrow \infty}e^{-rt}Y^y_t]  \leq  {\liminf}_{t \rightarrow \infty}\mathbb{E}[e^{-rt} Y^y_t]=0,$$
hence $\lim_{t \rightarrow \infty}e^{-rt}Y^y_t = 0$ $\P$-a.s.
\end{Remark}

In light of Assumption \ref{ass:ui} from now on we will adopt the convention
\begin{align}
\label{conventions}
e^{-r\tau}Y^y_{\tau}\mathds{1}_{\{\tau=\infty\}} := \lim_{t\rightarrow \infty}e^{-rt}Y^y_t =0, \quad a.s.
\end{align}
Also we set
\begin{align}
\label{conventions1}
e^{-r\tau}|f(X^x_{\tau},Y^y_{\tau})|\mathds{1}_{\{\tau=\infty\}}:=\limsup_{t\rightarrow \infty}e^{-rt}|f(X^x_{t},Y^y_{t})|, \quad a.s.,
\end{align}
for any Borel-measurable function $f$.

The next lemma will be useful in what follows.
\begin{Lemma}\label{ibp}
Under Assumptions \ref{ass:D2}, \ref{Ass:r} and \ref{ass:ui} it holds
\begin{align}
\label{integroparti}
\mathbb{E}[e^{-r \tau}Y^y_{\tau}]= y +\mathbb{E}\left[ \int_0^\tau e^{-r t}\big(\mu_2(Y^y_t)-r Y^y_t\big)dt\right],\quad\text{for $\tau\in \mathcal{T}$}.
\end{align}
\end{Lemma}
\textbf{Proof.}
The result holds for bounded stopping times $\tau_n:=\tau\wedge n$, with $\tau\in\mathcal{T}$ and $n\in\mathbb{N}$, by applying It\^o's formula, by noting that the resulting local martingale term is actually a true martingale by Assumptions \ref{ass:D2} and \ref{Ass:r} and by taking expectations. Then letting $n\rightarrow\infty$ and using Assumptions \ref{ass:D2}, \ref{ass:ui} and dominated convergence one finds \eqref{integroparti}.
\ep

In the rest of this section we aim at characterizing $v$ of \eqref{P:OS}.
\begin{Proposition}
\label{prop:v}
Under Assumptions \ref{ass:D2}, \ref{ass:cost}, \ref{Ass:r} and \ref{ass:ui} the following hold:
\begin{enumerate}
\item $v$ is such that
\begin{align}
\label{eq:subpon}
-y \leq v(x,y;z)\le C(z)(1+|x|^{\beta}+|y|), \qquad \forall(x,y)\in Q,
\end{align}
for a constant $C(z)>0$ depending on $z$.
\item $v(\,\cdot\,,y;z)$ is nonincreasing  for every $y\in \mathcal{I}_2$.
\item $v(x,\,\cdot\,;z)$ is nonincreasing  for every $x\in \mathcal{I}_1$.
\end{enumerate}
\end{Proposition}
\textbf{Proof.}
1. The lower bound follows by taking $\tau=0$ in \eqref{P:OS}. Assumptions \ref{ass:D2}, \ref{ass:cost}-$(iii)$, \ref{Ass:r}, \ref{ass:ui} and Lemma \ref{ibp} guarantee the upper bound.

2. The fact that $x\mapsto c_z(x,z)$ is nonincreasing (cf.\ Assumption \ref{ass:cost}-$(ii)$) and \eqref{comparison} imply
    \begin{align*}
    \hspace{-8pt}v(x_2,y;z)-v(x_1,y;z)\le \sup_{\tau\in\mathcal{T}}\E\Big[\int^\tau_0{e^{-rt}\big(c_z(X^{x_2}_t,z)-c_z(X^{x_1}_t,z)\big)dt}\Big]\le0,\quad\text{for $x_2>x_1$.}
    \end{align*}

3. It follows from \eqref{comparisonY} and arguments as in point 2.\ep

\begin{Proposition}
\label{vcont}
Under Assumptions \ref{ass:D2}, \ref{ass:cost}, \ref{Ass:r} and \ref{ass:ui} the value function $v(\,\cdot\,;z)$ of the optimal stopping problem \eqref{P:OS} is continuous on $Q$.
\end{Proposition}
\textbf{Proof.} Fix $z\in\mathbb{R}^+$ and let $\{(x_n,y_n), \,n\in\mathbb{N}\} \subset  Q$ be a sequence converging to $(x,y)\in Q$. Take $\varepsilon>0$ and let $\tau^\varepsilon:=\tau^\varepsilon(x,y;z)$ be an $\varepsilon$-optimal stopping time for the optimal stopping problem with value function $v(x,y;z)$. Then we have
\begin{align}
\label{pasd}
v(x,y;z)-v(x_n,y_n;z) \leq \varepsilon+\mathbb{E}\bigg[\int_0^{\tau^\varepsilon} e^{-rt} \big(c_z(X_t^x,z)-c_z(X^{x_n}_t,z)\big) dt-e^{-r\tau^\eps}(Y^y_{\tau^\varepsilon}- Y^{y_n}_{\tau^\varepsilon})\bigg].
\end{align}
Taking into account \eqref{contdata} and \eqref{contdataY}, Assumptions \ref{ass:cost}, \ref{Ass:r} and \ref{ass:ui}, we can apply dominated convergence (in its weak version requiring only convergence in measure; see, e.g., \cite[Ch.\ 2, Th.\ 2.8.5]{Bogachev}) to the right hand side of the inequality above and get
\begin{align}\label{eq:cont0}
\liminf_{n\rightarrow\infty}v(x_n,y_n;z)\geq v(x,y;z)-\varepsilon.
\end{align}

Similarly, taking $\varepsilon$-optimal stopping times $\tau_{n}^\varepsilon:=\tau^\varepsilon(x_n,y_n;z)$ for the optimal stopping pro\-blem with value function $v(x_n,y_n;z)$, and using Lemma \ref{ibp} we get
\begin{align}\label{limit2}
v(x_n,y_n;z)-v(x,y;z) \leq& \,\varepsilon+\mathbb{E}\bigg[\int_0^{\tau_{n}^\varepsilon} e^{-rt} \big(c_z(X_t^{x_n},z)-c_z(X^{x}_t,z)\big) dt-e^{-r\tau_{n}^\varepsilon}\big(Y^{y_n}_{\tau_{n}^\varepsilon}- Y^{y}_{\tau_{n}^\varepsilon}\big)\bigg]\nonumber\\
=& \,\varepsilon + \mathbb{E}\bigg[\int_0^{\tau_{n}^\varepsilon} e^{-rt} \big(c_z(X_t^{x_n},z)-c_z(X^{x}_t,z)\big)dt \bigg] - (y_n - y) \nonumber \\
 & + \mathbb{E}\bigg[\int_0^{\tau_{n}^\varepsilon} { e^{-r t} \big[r\big(Y^{y_n}_t - Y^{y}_t \big) - \big(\mu_2(Y^{y_n}_t) - \mu_2(Y^{y}_t)\big)\big]dt }\bigg] \\
 \leq & \,\varepsilon + \mathbb{E}\bigg[\int_0^{\infty} e^{-rt} \big|c_z(X_t^{x_n},z)-c_z(X^{x}_t,z)\big| dt \bigg]  + |y - y_n| \nonumber \\
&  + C\,\mathbb{E}\bigg[\int_0^{\infty} e^{-r t} \big|Y^{y_n}_t - Y^{y}_t \big| dt \bigg], \nonumber
\end{align}
for some $C>0$ and where we have used Lipschitz continuity of $\mu_2$ (cf.\ Assumption \ref{ass:D2}) in the last step.
Recalling now \eqref{contdata} and \eqref{contdataY}, \eqref{stimaKrylovXY}, Assumptions \ref{ass:cost} and \ref{Ass:r}, we can apply again dominated convergence in its weak version (cf.~\cite[Ch.\ 2, Th.\ 2.8.5]{Bogachev}) to the right hand side of the inequality above to obtain
\begin{align}\label{eq:cont1}
\limsup_{n\rightarrow\infty}v(x_n,y_n;z)\leq v(x,y;z)+\varepsilon.
\end{align}
Now \eqref{eq:cont0} and \eqref{eq:cont1} imply continuity of $v(\,\cdot\,,\,\cdot\,;z)$ by arbitrariness of $\varepsilon>0$.
\ep
\begin{Remark}
\label{rem:cont}
Arguments similar to those used in the proof of Proposition \ref{vcont} above may also be employed to show that $(x,y,z)\mapsto v(x,y;z)$ is continuous in $\mathcal{O}$.
\end{Remark}

We now provide a probabilistic representation of $v$, which we will later use to characterize the optimal boundary. For that, we first define the continuation and stopping regions of problem \eqref{P:OS} as
\beq
\label{czaz}
\mathcal{C}_z:=\{(x,y)\in Q \ | \ v(x,y;z)>-y\},\quad \mathcal{A}_z:=\{(x,y)\in Q \ | \ v(x,y;z)=-y\}.
\enq
We also recall that since $v(\,\cdot\,;z)$ is continuous standard optimal stopping theory (see, e.g., \cite{PeskShir}) guarantees that the stopping time
\begin{eqnarray}
\label{superarmonic2}
\tau^*=\tau^*(x,y;z):= \inf\,\{\,t\geq 0 \ | \ (X^x_t,Y^y_t) \in \Ac_z \,\}
\end{eqnarray}
is optimal for problem \eqref{P:OS}, whenever it is $\P$-a.s.~finite.
Moreover, we make the following
\begin{Assumption}
\label{ass:laws}
For every $(x,y)\in\Ic_1\times \Ic_2$ and $t > 0$ the laws of $X^x_t$ and $Y^y_t$ have densities $p_1(t,x,\,\cdot\,)$ and $p_2(t,y,\,\cdot\,)$, respectively. Moreover
\begin{itemize}
\item[1)] $(t,\zeta,\xi)\mapsto p_i(t,\zeta,\xi)$ is continuous on $(0,\infty)\times\Ic_i\times\Ic_i$, $i=1,2$;
\item[2)] For any compact set $\Kc\subset\Ic_1\times\Ic_2$ there exists $q>1$ (possibly depending on $\Kc$) such that 
\begin{align*}
\int_0^\infty e^{-rt} \left(\int_\Kc \big|p_1(t,x,\xi)p_2(t,y,\zeta)\big|^qd\xi d\zeta\right)^{\tfrac{1}{q}}dt<+\infty,\quad\text{for all $(x,y)\in\Kc$.}
\end{align*}
\end{itemize}
\end{Assumption}
\begin{Remark}
\label{rem:densities}
Assumption \ref{ass:laws} is clearly satisfied in the benchmark case of $X$ and $Y$ given by two independent geometric Brownian motions. The literature on the existence and smoothness of densities for the probability laws of solutions of SDEs driven by Brownian motion is huge and it mainly relies on PDEs' and Malliavin Calculus' techniques (see, e.g., \cite{Friedman} and \cite{Nualart} as classical references on the topic). In general, the existence of a density for the law of a one-dimensional diffusion is guaranteed under some very mild assumptions (see, e.g., the recent paper \cite{PrintemsFournier}). Sufficient conditions on our $(\mu_i,\sigma_i)$, $i=1,2$, to obtain Gaussian bounds for the transition densities and their first derivatives may be found for instance in \cite[Ch.~1, Th.~11]{Friedman}. One can also refer to, e.g., \cite{DeMarco} and references therein for more recent generalizations under weaker assumptions.
\end{Remark}

The proof of the next theorem is obtained in a number of steps, which we account for in the following technical subsection. Although these details are important, they are not necessary for the understanding of Sections \ref{sec:char} and \ref{sec:opt} and could be skipped at a first reading. 

\begin{Theorem}
\label{teo:formula}
Under Assumptions \ref{ass:D2}, \ref{ass:cost}, \ref{Ass:r}, \ref{ass:ui}, and \ref{ass:laws}, the following representation holds for every $(x,y)\in Q$:
\begin{align}
\label{v:rep}
v(x,y;z)=\mathbb{E}\bigg[\int_0^\infty e^{-rt}\left(c_z(X^x_t,z)\mathds{1}_{\{(X^x_t,Y^y_t)\in\mathcal{C}_z\}}-(rY^y_t-\mu_2(Y^y_t))\mathds{1}_{\{(X^x_t,Y^y_t)
\in\mathcal{A}_z\}}\right)\!dt\bigg].
\end{align}
\end{Theorem}   
Set
\beq\label{def:h}
H(x,y;z):=c_z(x,z)\mathds{1}_{\{(x,y)\in\Cc_z\}}-(ry-\mu_2(y))\mathds{1}_{\{(x,y)\in\Ac_z\}}
,\enq
so that \eqref{v:rep} can be written as
\beq\label{v:rep-short}
v(x,y;z)=\E\bigg[\int^\infty_0{e^{-r t}H(X^x_t,Y^y_t;z)dt}\bigg].
\enq

\noindent Due to \eqref{eq:subpon} and Assumption \ref{Ass:r}, the strong Markov property and standard arguments based on conditional expectations applied to the representation formula \eqref{v:rep-short} allow to verify that, for all $(x,y)\in Q$, 
\beq
\label{eq:mg-2}
e^{-r\tau}v(X^x_{\tau},Y^y_{\tau};z)+ \int^{\tau}_0 e^{-rs} H(X_s^x,Y_s^y;z)ds = \E\bigg[\int^\infty_0 e^{-rs}\,H(X^x_s,Y^y_s;z)\,ds\,\Big|\,\mathcal{F}_{\tau}\bigg],\,\tau\in\Tc,
\enq
and, in particular,
\beq
\label{eq:mg}
\bigg\{e^{-rt}v(X^x_t,Y^y_t;z)+ \int^t_0 e^{-rs} H(X_s^x,Y_s^y;z)ds,\, t \geq 0 \bigg\} \  \mbox{is an} \ (\Fc_t)\mbox{-martingale}.
\enq
Equation \eqref{eq:mg-2} also implies
\begin{align}\label{uiboundv}
\big|e^{-r\tau}v(X^x_\tau,Y^y_\tau;z)\big|\le\E\bigg[\int^\infty_0
{e^{-rt}\big|\,H(X^x_t,Y^y_t;z)\big|\,dt}\Big|\,\Fc_\tau\bigg]\,,\quad\tau\in\Tc,
\end{align}
and hence
the family $\big\{e^{-r\tau}v(X^x_\tau,Y^y_\tau;z)\,,\,\tau\in\Tc\big\}$ is uniformly integrable.

\subsection{Probabilistic representation of $v$: details}
Since the state space $Q=\mathcal{I}_1\times\mathcal{I}_2$ of the diffusion $\{(X^x_t,Y^y_t), t\geq 0\}$ may be unbounded, it is convenient for studying the variational inequality associated to our optimal stopping problem, to approximate problem \eqref{P:OS} by a sequence of problems on bounded domains. Let $\{Q_n, \, n \in  \mathbb{N}\}$ be a sequence of sets approximating $Q$, and we assume that
\beq\label{OOR}
\begin{cases}
Q_n\ \mbox{is open, bounded and  connected for every} \ n\in\N, \\
\partial{Q}_n\in C^{2+\alpha_n}\,\,\,\text{for some $\alpha_n>0$,}\\ 
{Q}_n\subset {Q}_{n+1} \ \mbox{for every}\   n\in\N,\\
 \lim_{n\rightarrow\infty} {Q}_n:=\bigcup_{n\geq 0}{Q}_n={Q}.
\end{cases}
\enq
Clearly, it is always possible to find such a sequence of sets.
The optimal stopping problem \eqref{P:OS} is then localized as follows.
Given $n\in\N$, define the stopping time
\beq
\label{taun}
\sigma_n=\sigma_n(x,y;z):=\inf\{t\geq 0\ | \ (X^x_t,Y^y_t) \notin {Q}_n\}
\enq
and notice that $\sigma_\infty=\sigma_\infty(x,y;z):=\inf\{t\geq 0\ | \ (X^x_t,Y^y_t)\notin {Q}\}=\infty$ a.s., since we are assuming that the boundaries of $X^x$ are non-exit and those of $Y^y$ are natural, hence non attainable. Moreover, from the last of \eqref{OOR} we obtain
\beq\label{tauRtau}
\sigma_n\uparrow \sigma_\infty=\infty \quad \mbox{$\P$-a.s., \ as } n\rightarrow \infty.
\enq
With $\sigma_n$ as in \eqref{taun}, we can define the approximating optimal stopping problem
\beq\label{P:OSn}
v_n(x,y;z):=\sup_{\tau\in\mathcal{T}} \mathbb{E}\bigg[\int_{0}^{\sigma_n\wedge\tau}e^{-rt} c_z(X^x_t, z)dt-e^{-r(\sigma_n\wedge\tau)}Y^y_{\sigma_n\wedge\tau} \bigg], \ \ \ (x,y)\in {Q},
\enq
and prove the following

\begin{Proposition}
\label{prop:vR}
Let Assumptions \ref{ass:D2}, \ref{ass:cost}, \ref{Ass:r} and \ref{ass:ui} hold. Then
\begin{enumerate}
\item $v_n(\,\cdot\,;z)\le v_{n+1}(\,\cdot\,;z)\leq v(\,\cdot\,;z)$ on ${Q}$ for all $n\in\mathbb{N}$.
\item
 $v_n(x,y;z)=-y$ for $(x,y)\in Q\,\backslash\, Q_n$ and all $n\in\mathbb{N}$ (in particular for every $(x,y)\in \partial Q_n$, since $Q_n$ is open).
\item $v_n(x,y;z)\uparrow v(x,y;z)$ as $n\rightarrow\infty$ for every $(x,y)\in {Q}$.
\item If $\{v_n(\,\cdot\,;z), \,n\in\mathbb{N}\}\subset C^0({Q})$, then $v_n(\,\cdot\,;z)$ converges to $v(\,\cdot\,;z)$ uniformly on all compact subsets $\mathcal{K} \subset {Q}$.
\end{enumerate}
\end{Proposition}
\noindent\textbf{Proof.}
1. It follows from \eqref{tauRtau} and by comparison of \eqref{P:OSn} with \eqref{P:OS}.

2. This claim follows from the definition of $\sigma_n$ and of $v_n$ (see \eqref{taun} and \eqref{P:OSn}, respectively).

3. For fixed $(x,y)\in Q$ denote by $\tau^\varepsilon:=\tau^\varepsilon(x,y;z)$ an $\varepsilon$-optimal stopping time of $v(x,y;z)$, then
\begin{eqnarray*}
0 \hspace{-0.25cm} & \leq & \hspace{-0.25cm} v(x,y;z)-v_n(x,y;z)\\
\hspace{-0.25cm} & \leq & \hspace{-0.25cm} \mathbb{E}\bigg[\int_{\sigma_n\wedge \tau^\varepsilon} ^{\tau^\varepsilon}e^{-rt}c_z(X^x_t, z)dt
 - \left(e^{-r \tau^\varepsilon}Y^y_{\tau^\varepsilon}- e^{-r\sigma_n}Y^y_{\sigma_n}\right)\mathds{1}_{\{\sigma_n<\tau^\varepsilon\}}\bigg] + \varepsilon,
\end{eqnarray*}
where the first inequality is due to $1$\ above. Now, the sequence of random variables $\{Z_n, n\in\N\}$ defined by
$$Z_n:= \int_{\sigma_n\wedge\tau^\varepsilon} ^{\tau^\varepsilon}e^{-rt}c_z(X^x_t, z)dt-\left(e^{-r \tau^\varepsilon}Y^y_{\tau^\varepsilon}- e^{-r\sigma_n}Y^y_{\sigma_n}\right)\mathds{1}_{\{\sigma_n< \tau^\varepsilon\}}$$
is uniformly integrable due to Assumptions \ref{ass:cost}, \ref{Ass:r} and \ref{ass:ui}, and $\lim_{n \rightarrow \infty} Z_n= 0$ $\P$-a.s., by Remark \ref{Rem:ass:ui}-(\textbf{2}) and \eqref{tauRtau}. Then $3$\ follows from Vitali's convergence theorem and arbitrariness of $\varepsilon$.

4.
Since $v(\,\cdot\,;z)\in C^0({Q})$, the claim follows from $1$~and $3$~above and by Dini's Lemma.
\ep

\smallskip
Fix $n\in\N$ and $z\in\R^+$, and define the continuation and stopping regions of our approximating optimal stopping problem \eqref{P:OSn} respectively by
\begin{align}\label{def:CnAn}
\mathcal{C}^n_z:=\{(x,y)\in Q\ | \ v_n(x,y;z)>-y\}, \quad \mathcal{A}^n_z:=\{(x,y)\in Q \ | \ v_n(x,y;z)=-y\}.
\end{align}

Denote by $\mathbb{L}$ the second order elliptic differential operator associated to the two-dimensional diffusion $\{(X_t,Y_t), t\ge0\}$. Since $X$ and $Y$ are independent then $\mathbb{L}:=\mathbb{L}_X+\mathbb{L}_Y$, with
\begin{align*}
\vspace{+8pt}
&(\mathbb{L}_Xf)\,(x,y):=\frac{1}{2}(\sigma_1)^2(x)\frac{\partial^2}{\partial x^2}f(x,y)+\mu_1(x)\frac{\partial}{\partial x}f(x,y),\\
&(\mathbb{L}_Yf)\,(x,y):=\frac{1}{2}(\sigma_2)^2(y)\frac{\partial^2}{\partial y^2}f(x,y)+\mu_2(y)\frac{\partial}{\partial y}f(x,y),
\end{align*}
for $f\in C^2_b(\overline{Q})$. From standard arguments we can formally associate the function $v_n(\,\cdot\,,\,\cdot\,;z)|_{Q_n}$ of \eqref{P:OSn} to the variational inequality (parametrized in $z$)
\begin{align}
\label{VIn}
\max\Big\{\big(\mathbb{L}-r\big)u(x,y;z)+c_z(x,z), -u(x,y;z)-y\Big\}=0, \ \ \ (x,y)\in{Q}_n,
\end{align}
with boundary condition
\begin{align}\label{BC}
u(x,y;z)=-y, \ \ \ \ (x,y)\in\partial Q_n.
\end{align}

The next result is standard and its proof is given in the Appendix for the sake of completeness.
\begin{Proposition}
\label{teo:ex}
Under Assumptions \ref{ass:D2}, \ref{ass:cost}, \ref{Ass:r} and \ref{ass:ui}, for each $n\in\N$ and $z\in\R^+$ $v_n(\cdot\,;z)\in W^{2,p}(Q_n)$ for all $1\leq p<\infty$, and uniquely solves \eqref{VIn} a.e.\ in $Q_n$ with the boundary condition \eqref{BC}. Moreover, the stopping time
\beq\label{opt-st-n}
\tau^*_{n}(x,y;z):=\inf\big\{t\ge0\,|\,(X^x_t,Y^y_t)\notin \Cc^n_z \big\},
\enq
with $\Cc^n_z$ as in \eqref{def:CnAn}, is optimal for problem \eqref{P:OSn}.
\end{Proposition}

\begin{Remark}
\label{Sobolev}
Note that, by well known Sobolev's inclusions (see for instance \cite[Ch.\ 9, Cor.\ 9.15]{Br}), the space $W^{2,p}(Q_n)$ with $p \in(2,\infty)$ can be continuously embedded into  $C^1(\overline{Q}_n)$. Hence,  the boundary condition \eqref{BC} is well-posed for functions in the class $W^{2,p}(Q_n)$, $p\in(2,\infty)$. In the following we shall always refer to the unique $C^1$ representative of elements of $W^{2,p}(Q_n)$. 
\end{Remark}

\begin{Proposition}
For every $(x,y)\in Q$ the following representation holds
\begin{align}\label{vn:rep}
v_n(x,y;z)  = \mathbb{E}\bigg[\hspace{-1pt}\int_0^{\sigma_n}\hspace{-4pt}e^{-rt}\hspace{-2pt}\left(\hspace{-1pt}c_z(X^x_t,z)
\mathds{1}_{\{(X^x_t,Y^y_t)\in\mathcal{C}^n_z\}}\hspace{-2pt}
-\hspace{-2pt}(rY^y_t\hspace{-1pt}-\hspace{-1pt}\mu_2(Y^y_t))\mathds{1}_{\{(X^x_t,Y^y_t)\in\mathcal{A}^n_z\}}\hspace{-1pt}\right)dt\hspace{-1pt} -e^{-r\sigma_n}\hspace{-1pt} Y^y_{\sigma_n}
\bigg].
\end{align}
\end{Proposition}
\textbf{Proof.}
Since $v_n(\,\cdot\,;z) \in W^{2,p}(Q_n)$ and solve \eqref{VIn}--\eqref{BC} (cf.\ Proposition \ref{teo:ex}), a generalised It\^o's formula gives (see also \eqref{dyn02} and \eqref{zero} in the Appendix)
\begin{align}
\label{dyn021}
{v}_n(x,y;z)=\E\left[-e^{-r\sigma_n}
Y^y_{\sigma_n}-\int^{\sigma_n}_0{e^{-rt}(\mathbb{L}-r){v}_n(X^x_t,Y^y_t;z)\,dt}
\right].
\end{align}
It follows from Propositions \ref{teo:ex} that\footnote{There is a small technicality concerning this claim, which we account for in Lemma \ref{Lemmahessiano} of Appendix \ref{LemmaAppendix} for the interested reader.}  
\begin{align}\label{var01}
(\mathbb{L}-r)v_n(x,y;z)=c_z(x,z)\mathds{1}_{\{(x,y)\in\Cc^n_z\}}-(ry-\mu_2(y))\mathds{1}_{\{(x,y)\in\Ac^n_z\}},
\quad\text{for a.e.~$(x,y)\in Q_n$,}
\end{align}
and we have the claim by using \eqref{dyn021} and Assumption \ref{ass:laws} in \eqref{var01}.
\ep\\

We observe that since $v_n \leq v$ and $\{v_n, n\in\N\}$ is an increasing sequence then
\beq\label{cnan1}
\mathcal{C}^n_z\subset \mathcal{C}^{n+1}_z\subset \mathcal{C}_z, \quad \mathcal{A}^n_z\supset\mathcal{A}^{n+1}_z\supset\mathcal{A}_z, \ \ \ \forall n\in\N.
\enq
On the other hand, the pointwise convergence $v_n \uparrow v$ (cf.\ Proposition \ref{prop:vR}) implies that if $(x_0,y_0)\in\Cc_z$, then $v(x_0,y_0)+y_0\ge\varepsilon_0$ for some $\varepsilon_0>0$ and $v_n(x_0,y_0)+y_0\ge\varepsilon_0/2$ for all $n\ge n_0$ and suitable $n_0\in\N$. Hence we have
\beq\label{cnan}
\lim_{n\rightarrow\infty}\mathcal{C}^n_z:=\bigcup_{n\geq 0} \mathcal{C}^n_z=\mathcal{C}_z, \ \ \ \ \lim_{n\rightarrow\infty}\mathcal{A}^n_z:=\bigcap_{n\geq 0} \mathcal{A}_z^n=\mathcal{A}_z.
\enq
We can now prove Theorem \ref{teo:formula}. 

\textbf{Proof of Theorem \ref{teo:formula}.}
We study \eqref{vn:rep} in the limit as $n\uparrow \infty$. 
Observe that:

\emph{1.} The left-hand side of \eqref{vn:rep} converges pointwisely to $v(x,y;z)$ by Proposition \ref{prop:vR}-(3);

\emph{2.} $\{e^{-r\sigma_n}Y^y_{\sigma_n},\,n\in\mathbb{N}\}$ is a family of random variables uniformly integrable and converging a.s.\ to $0$, due to \eqref{tauRtau} and to Assumptions \ref{Ass:r} and \ref{ass:ui} (see also the discussion in Remark \ref{Rem:ass:ui}-(\textbf{2})). Hence $\lim_{n \rightarrow \infty}\mathbb{E}\left[e^{-r\sigma_n}Y^y_{\sigma_n}\right]= 0,$ by Vitali's convergence Theorem;

\emph{3.} From \eqref{cnan1}, one has
\begin{align}\label{conv01}
& \left|\,\mathbb{E}\bigg[\int_0^{\sigma_n} e^{-rt}c_z(X^x_t,z)\mathds{1}_{\{(X_t^x,Y^y_t)\in\mathcal{C}_n\}}dt- \int_0^\infty e^{-rt}c_z(X^x_t,z)\mathds{1}_{\{(X_t^x,Y^y_t)\in\mathcal{C}\}}dt\bigg]\right| \\
& \leq \mathbb{E}\bigg[\int_0^{\infty} e^{-rt}|c_z(X^x_t,z)|\mathds{1}_{\{(X_t^x,Y^y_t)\in\mathcal{C}\,\setminus\, \mathcal{C}_n\}}dt\bigg]+\mathbb{E}\bigg[\int_{\sigma_n}^\infty e^{-rt}|c_z(X^x_t,z)|\mathds{1}_{\{(X_t^x,Y^y_t)\in\mathcal{C}\}}dt\bigg].\nonumber
\end{align}
The first term in the right-hand side of \eqref{conv01} converges to zero as $n\rightarrow\infty$ by dominated convergence and \eqref{cnan} (cf.~Assumptions \ref{ass:cost}-(iii), \ref{Ass:r} and Remark \ref{rem:c}-(\textbf{3})).
Similarly, dominated convergence and \eqref{tauRtau} give
$$\lim_{n \rightarrow \infty}\mathbb{E}\bigg[\int_{\sigma_n}^\infty e^{-rt}|c_z(X^x_t,z)|\mathds{1}_{\{(X_t^x,Y^y_t)\in\mathcal{C}\}}dt\bigg] = 0.$$

\emph{4.} From \eqref{cnan} it follows that for a.e.\ $(t,\omega)\in\mathbb{R}^+\times \Omega$
\begin{align*}
\lim_{n \rightarrow \infty }\mathds{1}_{[0,\sigma_n]}(t) e^{-rt}\Big[rY^y_t-\mu_2(Y^y_t)\Big]\mathds{1}_{\{(X^x_t,Y^y_t)\in\mathcal{A}^n_z\}} = e^{-rt}\Big[rY^y_t-\mu_2(Y^y_t)\Big]\mathds{1}_{\{(X^x_t,Y^y_t)\in \mathcal{A}_z\}}.
\end{align*}
Moreover, due to  Lipschitz-continuity of $\mu_2$ (cf.\ Assumption \ref{ass:D2}),
\begin{eqnarray*}
\Big|e^{-rt}\Big[rY^y_t-\mu_2(Y^y_t)\Big]\mathds{1}_{\{(X^x_t,Y^y_t)\in\mathcal{A}^n_z\}}\Big|\  \leq \ e^{-rt}\Big|rY^y_t-\mu_2(Y^y_t)\Big| \  \leq \ e^{-rt} C_0(1 + Y^y_t),
\end{eqnarray*}
for some $C_0 > 0$ depending on $y$ and $r$.
The last expression of the inequality above is integrable in $\mathbb{R}^+\times \Omega$ by \eqref{stimaKrylovXY} and by Assumption \ref{Ass:r}. Hence dominated convergence and \eqref{tauRtau} yield
\begin{align*}
\lim_{n \rightarrow \infty}\mathbb{E}\bigg[\hspace{-2pt}\int_0^{\sigma_n}\hspace{-5pt} e^{-rt}\Big[rY^y_t-\mu_2(Y^y_t)\Big]\mathds{1}_{\{(X^x_t,Y^y_t)\in\mathcal{A}^n_z\}}dt\bigg]= \mathbb{E}\bigg[\hspace{-2pt}\int_0^\infty\hspace{-5pt} e^{-rt}\Big[rY^y_t-\mu_2(Y^y_t)\Big]\mathds{1}_{\{(X^x_t,Y^y_t)\in\mathcal{A}_z\}}dt\bigg].
\end{align*}
\medskip
Now taking $n\rightarrow \infty$ in \eqref{vn:rep} and using \emph{1-4} above, \eqref{v:rep} follows.
\ep \\


\section{Characterization of the Optimal Boundary}
\label{sec:char}

In this section we will provide a characterization of the optimal boundaries of the family of optimal stopping problems \eqref{P:OS}. For that we define
\begin{align}
\label{ystar}
y^*(x;z):=\inf \{y\in\mathcal{I}_2\ | \ v(x,y;z)>-y\}, \ \ \ (x,z)\in\mathcal{I}_1\times\mathbb{R}^+,
\end{align}
with the convention $\inf\emptyset = \overline{y}$. Notice that under this convention $y^*(\,\cdot\,;z)$ takes values in $\overline{\Ic}_2$.
We will show that under suitable conditions $y^*(\,\cdot\,;z)$ splits $\Ic_1\times\Ic_2$ into $\Cc_z$ and $\Ac_z$ (cf.\ \eqref{czaz}). Moreover, we will characterize $y^*(\,\cdot\,;z)$ as the unique continuous solution of a nonlinear integral equation of Fredholm type.
\begin{Remark}
A common way of obtaining integral equations of optimal stopping boundaries is by using the so-called local time space formula (cf.~\cite{Peskir}). In our case this would require to prove that the process $\{y^*(X^x_t;z),\,t\ge0\}$ is a semimartingale for each given $z \in \R^+$ (see \cite[Th.\ 2.1]{Peskir}). Proving the latter is extremely challenging. Here we obtain the same integral equation but following a different approach that builds on results of the previous section.
\end{Remark}
Throughout this section Assumptions \ref{ass:D2}, \ref{ass:cost}, \ref{Ass:r}, \ref{ass:ui}, and \ref{ass:laws} will be standing assumptions and we will not repeat them in the statements of the next results. 
We now make the following
\begin{Assumption}
\label{ass:mu21}
The map $y\mapsto ry-\mu_2(y)$ is increasing.
\end{Assumption}

\begin{Proposition}
\label{prop:monot}
Under Assumption \ref{ass:mu21} one has (cf.~\eqref{czaz})
\beq\label{ACz}
\mathcal{C}_z=\{(x,y)\in Q\ | \ y > y^*(x;z)\}, \ \ \ \mathcal{A}_z=\{(x,y)\in Q\ | \ y \leq y^*(x;z)\}.
\enq
\end{Proposition}
\textbf{Proof.}
It suffices to show that $y\mapsto v(x,y;z)+y$ is nondecreasing for each $x\in \mathcal{I}_1$, $z\in\R^+$. Set $\bar{u}:=v+y$, take $y_1$ and $y_2$ in $\Ic_2$ such that $y_2>y_1$ and set $\tau_1:=\inf\{t\ge0\ | \ (X^x_t,Y^{y_1}_t)\notin \Cc_z\}$, which is optimal for $v(x,y_1;z)$. From Lemma \ref{ibp}, the well known superharmonic characterization of $v$ 
and \eqref{comparisonY} we obtain
\begin{align}\label{ubar02}
\bar{u}(x,y_2;z)-\bar{u}(x,y_1;z)\ge&\ \,\E\bigg[e^{-r{\tau_1}}\big(\bar{u}(X^x_{{\tau_1}},Y^{y_2}_{{\tau_1}};z)-
\bar{u}(X^x_{{\tau_1}},Y^{y_1}
_{{\tau_1}};z)\big)\bigg]\nonumber\\
&\ +\E\bigg[\int^{{\tau_1}}_0{e^{-rt}\Big(r\big(Y^{y_2}_t-Y^{y_1}_t\big)-\big(\mu_2(Y^{y_2}_t)-\mu_2(Y^{y_1}_t)\big)\Big)dt}
\bigg]\\
\ge&\,\ \E\bigg[e^{-r{\tau_1}}\left(\bar{u}(X^x_{{\tau_1}},Y^{y_2}_{{\tau_1}};z)-
\bar{u}(X^x_{{\tau_1}},Y^{y_1}
_{{\tau_1}};z)\right)\bigg],\nonumber
\end{align}
where the last inequality follows by \eqref{comparisonY} and Assumption \ref{ass:mu21}. Note that the last expression in \eqref{ubar02} is well defined thanks to Assumption \ref{ass:ui} and \eqref{uiboundv}. Moreover, since $\bar{u}\geq 0$ it holds
\begin{align}\label{ubar04bis}
\E\bigg[e^{-r{\tau_1}}\left(\bar{u}(X^x_{{\tau_1}},Y^{y_2}_{{\tau_1}};z)-
\bar{u}(X^x_{{\tau_1}},Y^{y_1}
_{{\tau_1}};z)\right)\bigg]\ge -\E\bigg[e^{-r{\tau_1}}
\bar{u}(X^x_{{\tau_1}},Y^{y_1}
_{{\tau_1}},z)\bigg].
\end{align}
By Assumption \ref{Ass:r}, Proposition \ref{prop:v}-(1) and since $\mathds{1}_{\{\tau_1\le n\}}e^{-r{\tau_1}}\bar{u}(X^x_{\tau_1},Y^{y_1}_{\tau_1};z)=0$ $\P$-a.s., Fatou's Lemma gives
\begin{align}\label{ubar04}
\E\bigg[e^{-r{\tau_1}}\bar{u}(X^x_{{\tau_1}},Y^{y_1}
_{{\tau_1}};z)\bigg]=&\ \E\bigg[\liminf_{n \rightarrow \infty} e^{-r({\tau_1}\wedge n)}\bar{u}(X^x_{{\tau_1}\wedge n},Y^{y_1}
_{{\tau_1}\wedge n};z)\bigg]\nonumber\\
\le& \ \liminf_{n\rightarrow\infty}\E\left[e^{-r n}\bar{u}(X^x_{ n},Y^{y_1}
_{ n};z)\mathds{1}_{\{{\tau_1}>n\}}\right]=0
\end{align}

Now \eqref{ubar02}, \eqref{ubar04bis}, and \eqref{ubar04} imply that $y\mapsto\bar{u}(x,y;z)$ is increasing,  therefore \eqref{ACz} holds.

\ep

Notice that \eqref{v:rep} and \eqref{ACz} imply
\beq\label{v:rep1}
v(x,y;z)=\mathbb{E}\bigg[\int_0^\infty e^{-rt}\left(c_z(X^x_t,z)\mathds{1}_{\{Y^y_t>y^*(X^x_t;z)\}}-(rY^y_t-\mu_2(Y^y_t))\mathds{1}_{\{Y^y_t\leq y^*(X^x_t;z)\}}\right)\!dt\bigg].
\enq
Under Assumption \ref{ass:laws}, \eqref{v:rep1} can also be expressed in a purely analytical way as
\begin{align}
\label{eq:intanv}
v(x,y;z) = &\int_0^\infty e^{-rt} \bigg[\int^{\overline{x}}_{\underline{x}}p_1(t,x,\xi)c_z(\xi,z)\bigg(\int^{\overline{y}}_{y^*(\xi;z)} \ p_2(t,y,\eta) d\eta\bigg)d\xi\bigg]dt\\
&  -\,\int_0^\infty e^{-rt} \bigg[\int^{\overline{x}}_{\underline{x}} p_1(t,x,\xi)\bigg(\int^{y^*(\xi;z)}_{\underline{y}} (r\eta-\mu_2(\eta )) p_2(t,y,\eta)d\eta\bigg)d\xi\bigg]dt,\nonumber
\end{align}
for any $(x,y,z) \in \mathcal{O}$.

\begin{Proposition}\label{prop:y*}
Under  Assumption \ref{ass:mu21} one has
\begin{enumerate}
\item the function $y^*(\,\cdot\,;z)$ is nondecreasing and right-continuous for any $z\in \mathbb{R}^+$;
\item the function $y^*(x;\,\cdot\,)$ is nonincreasing and left-continuous for any $x\in \mathcal{I}_1$;
\end{enumerate}
\end{Proposition}
\textbf{Proof.}
Claims 1 and 2 follow by adapting arguments from the proof of \cite[Prop.\ 2.2]{Jacka} and by using our Proposition \ref{prop:v}-(2)-(3), and Proposition \ref{vcont}.\ep

\medskip
It follows from Propositions \ref{prop:monot} and \ref{prop:y*}-(1) that the regions $\mathcal{C}_z$ and $\mathcal{A}_z$ are connected for every $z \in \mathbb{R}^+$, and the optimal stopping time $\tau^*(x,y;z)$ defined in  \eqref{superarmonic2} can be written as
\begin{align}
\label{deftauhat}
\tau^*(x,y;z)=\inf\big\{t\ge0\,|\,Y^y_t\le y^*(X^x_t;z)\big\}.
\end{align}

Thanks to the representation \eqref{v:rep1} or \eqref{eq:intanv}, under the following further assumptions we can prove the $C^1$-regularity of the function $v$.
\begin{Assumption}
\label{ass:densitiesderivative}
The functions $p_1(t,\cdot,\xi)$ and $p_2(t,\cdot,\eta)$ are differentiable for each $(t,\xi)\in (0,\infty) \times \mathcal{I}_1$ and each $(t,\eta)\in (0,\infty) \times \mathcal{I}_2$, respectively.
Moreover, denoting by $p^\prime_i$, $i=1,2$ the partial derivative of $p_i$ with respect to the second variable, it holds
\begin{itemize}
\item[1)] $x\mapsto p^\prime_1(t,x,\xi)$ is continuous in $\Ic_1$ for all $(t,\xi)\in (0,\infty) \times\Ic_1$ and, for any $(x,y,z)\in\Oc$, there exists $\delta>0$ such that $\sup_{\zeta\in[x-\delta,x+\delta]}\big|p^\prime_1(t,\zeta,\xi)\big|\le\psi_1(t,\xi;\delta)$ for some $\psi_1$ such that    
\begin{align}
\int_0^\infty e^{-rt}\left(\int_Q\psi_1(t,\xi;\delta)p_2(t,y,\eta)\big|c_z(\xi,z)+\eta\big|d\xi\,d\eta\right)dt<+\infty
\end{align}
\item[2)] $y\mapsto p^\prime_2(t,y,\eta)$ is continuous in $\Ic_2$ for all $(t,\eta)\in (0,\infty) \times\Ic_2$ and, for any $(x,y,z)\in\Oc$, there exists $\delta>0$ such that $\sup_{\zeta\in[y-\delta,y+\delta]}\big|p^\prime_2(t,\zeta,\eta)\big|\le\psi_2(t,\eta;\delta)$ for some $\psi_2$ such that
\begin{align}
\int_0^\infty e^{-rt}\left(\int_Q\psi_2(t,\eta;\delta)p_1(t,x,\xi)\big|c_z(\xi,z)+\eta\big|d\xi\,d\eta\right)dt<+\infty
\end{align}
    \end{itemize}
\end{Assumption}

\begin{Proposition}
\label{C1}
Under Assumptions \ref{ass:mu21} and \ref{ass:densitiesderivative}, one has $v(\,\cdot\,;z)\in C^1(Q)$ for every $z\in\mathbb{R}^+$.
\end{Proposition}
\textbf{Proof.}
The proof follows by \eqref{eq:intanv},  Assumption \ref{ass:densitiesderivative}, and standard dominated convergence arguments.
\ep
\vspace{+4pt}

Proposition \ref{C1} above states in particular the so-called  \emph{smooth-fit} condition across the free-boundary, i.e.\ the continuity of $v_x(\,\cdot\,;z)$ and $v_y(\,\cdot\,;z)$ at $\partial\mathcal{A}_z$.
With the aim of characterizing the boundary $y^*(\,\cdot\,;z)$ as unique \emph{continuous} solution of a (parametric) integral equation we make the following additional
\begin{Assumption}
\label{ass:mu2}
The drift coefficient $\mu_2$ is continuously differentiable in $\Ic_2$ and $ \frac{\partial\mu_2}{\partial y}< r$. Moreover, $\mu_2,\sigma_2\in C^{1+\delta}(\mathcal{I}_2)$, for some $\delta > 0$.
\end{Assumption}

\begin{Proposition}
\label{continuousboundary}
Under Assumptions \ref{ass:mu21}, \ref{ass:densitiesderivative} and \ref{ass:mu2}, the function $y^*(\,\cdot\,;z):\mathcal{I}_1\rightarrow \overline{\mathcal{I}}_2$ is continuous.
\end{Proposition}
\textbf{Proof.}
We know that the function $y^*(\,\cdot\,;z)$ is nondecreasing and right-continuous by Proposition \ref{prop:y*}-(1). Hence it suffices to show that it is also left-continuous. Borrowing arguments from \cite{DeA2015}, we argue by contradiction and we assume that there exists $x_0\in \mathcal{I}_1$ such that $y^*(x_0-;z):=\lim_{x\uparrow x_0} y^*(x;z)<y^*(x_0;z)$. Then, there also exist $y_0\in\mathcal{I}_2$ and $\varepsilon>0$ such that
$$
\Sigma_z:= (x_0-\varepsilon,x_0)\times (y_0-\varepsilon,y_0+\varepsilon)\subset \mathcal{C}_z, \ \ \ \  \{x_0\}\times (y_0-\varepsilon,y_0+\varepsilon)\subset \mathcal{A}_z.
$$
Notice that, by standard arguments on free-boundary problems and optimal stopping (cf.\ for instance \cite[Ch.~3, Sec.~7]{PeskShir} discussion at p.\ 131 together with PDE result \cite[Ch.\ 6, Sec.\ 3, Thm.6.13]{GT}), one has that $v(\,\cdot\,;z)\in C^2(\Cc_z)$ and solves
\begin{align}\label{eq:pde}
\frac{1}{2}\sigma_1^2(x)v_{xx}(x,y;z)=-\mu_1(x)v_x(x,y;z)-(\mathbb{L}_Y-r)v(x,y;z)-c_z(x,z), \ \ \ \ \ (x,y)\in\mathcal{C}_z.
\end{align}
On the other hand, since $\mu_2,\sigma_2\in C^{1+\delta}(\mathcal{I}_2)$, regularity results on uniformly elliptic partial differential equations (cf.~for instance \cite[Ch.\ 6, Th.\ 6.17]{GT}) imply that one actually has $v_y(\,\cdot\,;z)\in C^{2+\delta}(\mathcal{C}_z)$. Hence we can differentiate \eqref{eq:pde} with respect to $y$ to find
\begin{align}
\label{ppl}
\frac{1}{2}\sigma_1^2(x)(v_{y})_{xx}(x,y;z)=-\mu_1(x)(v_{y})_x(x,y;z)-(\mathcal{R}-r)v_y(x,y;z), \ \ \ \ \ (x,y)\in\mathcal{C}_z,
\end{align}
where
$$(\mathcal{R}f)(x,y):=\frac{1}{2}\sigma^2_2(y)f_{yy}(x,y)+\Big[\frac{\partial \sigma^2_2}{\partial y}(y)+\mu_2(y)\Big] f_{y}(x,y)+\frac{\partial\mu_2}{\partial y}(y)f(x,y), \ \ \ \ \ f\in C^2_b(Q).$$
Take now $y_1,y_2\in (y_0-\varepsilon,y_0+\varepsilon)$ with $y_1<y_2$ and set
\beq
\label{defFphi}
F_\phi(x; y_1,y_2,z):=-\int_{y_1}^{y_2}v_{xx}(x,y;z)\phi'(y)dy,\ \ \ \ x\in(x_0-\varepsilon,x_0),
\enq
where $\phi$ is real-valued, arbitrarily chosen and such that
$$\phi\in C_c^\infty(y_1,y_2),\ \ \ \ \phi\geq 0, \ \ \ \ \int_{y_1}^{y_2}\phi(y)dy>0.$$
From now on we will write $F_\phi(x)$ instead of $F_\phi(x; y_1,y_2,z)$ to simplify the notation. Multiply both sides of \eqref{ppl} by $2\phi(y)/\sigma^2_1(x)$ and integrate by parts with respect to $y\in(y_1,y_2)$; it follows
\beq\label{pkj}
F_\phi(x) \hspace{-0.25cm} & =& \hspace{-0.25cm} -\int_{y_1}^{y_2}\frac{1}{\sigma^2_1(x)}\Big[\mu_1(x)v_{xy}(x,y;z)+(\mathcal{R}-r)v_y(x,y;z)\Big]\phi(y)dy \\
\hspace{-0.25cm} &=& \hspace{-0.25cm} \frac{\mu_1(x)}{\sigma^2_1(x)}\int_{y_1}^{y_2}v_x(x,y;z)\phi'(y)dy+\frac{1}{\sigma_1^2(x)}\int_{y_1}^{y_2} v(x,y;z)\frac{\partial}{\partial y}(\mathcal{R}-r)^*\phi(y)dy,\nonumber
\enq
for every $x\in(x_0-\varepsilon,x_0)$, with $(\mathcal{R}-r)^*$ denoting the adjoint of $(\mathcal{R}-r)$. Now, recalling Proposition \ref{C1} and the definition of $\mathcal{C}_z$ and $\mathcal{A}_z$ one also has
\beq\label{sf}
\begin{cases}
v(x_0,y;z)= -y, \ \ \ \forall y\in [y_1,y_2],\\
v_x(x_0,y;z)= 0, \ \ \ \ \forall y\in [y_1,y_2],\\
v_y(x_0,y;z)= -1, \ \ \forall y\in [y_1,y_2].
\end{cases}
\enq
Thus, taking limits in \eqref{pkj}, one obtains
\beq\label{lim:F}
\lim_{x\uparrow x_0} F_\phi(x) \hspace{-0.25cm} &=  & \hspace{-0.25cm} -\frac{1}{\sigma_1^2(x_0)}\int_{y_1}^{y_2}y\frac{\partial}{\partial y}(\mathcal{R}-r)^*\phi(y)dy = \frac{1}{\sigma_1^2(x_0)}\int_{y_1}^{y_2} [(\mathcal{R}-r)1]\phi(y)dy \nonumber \\
\hspace{-0.25cm} & = & \hspace{-0.25cm} \frac{1}{\sigma_1^2(x_0)}\int_{y_1}^{y_2} \Big(\frac{\partial}{\partial y}\mu_2(y)-r\Big) \phi(y)dy \ <\ 0,
\enq
where the last inequality follows from  Assumption \ref{ass:mu2}. Since $F_\phi$ is clearly continuous in $(x_0-\varepsilon,x_0)$, we see from \eqref{lim:F} that it must be $F_\phi< 0$ in a left neighborhood of $x_0$ and, without any loss of generality, we  assume that $F_\phi< 0$ in $(x_0-\varepsilon,x_0)$. Recalling \eqref{defFphi}, we have for each $\delta\in(0,\varepsilon)$
\beqs
0 \hspace{-0.25cm} &  >  & \hspace{-0.25cm} \int_{x_0-\delta}^{x_0} F_{\phi}(x)dx = -\int_{x_0-\delta}^{x_0} \int_{y_1}^{y_2} v_{xx}(x,y;z)\phi'(y)dy\,dx\\
\hspace{-0.25cm} &  =  & \hspace{-0.25cm} -\int_{y_1}^{y_2} [v_{x}(x_0,y;z)-v_x(x_0-\delta,y;z)]\phi'(y)dy \\
\hspace{-0.25cm} &  = & \hspace{-0.25cm} \int_{y_1}^{y_2} v_x(x_0-\delta,y;z)\phi'(y)dy = -\int_{y_1}^{y_2}v_{xy}(x_0-\delta,y;z)\phi(y)dy,
\enqs
by \eqref{sf} and Fubini-Tonelli's theorem. This implies that $v_{xy}(\,\cdot\,;z)>0$ in $\Sigma_z$ by arbitrariness of $\phi$ and $\delta$ and hence the function $x\mapsto v_{y}(x,y;z)$ is strictly increasing in $(x_0-\varepsilon,x_0)$ for any $y\in[y_1,y_2]$. It then follows from the last of \eqref{sf}
\beq\label{vy<}
v_y(\,\cdot\,;z)<-1 \ \  \mbox{ in} \ \Sigma_z\subset \mathcal{C}_z.
\enq
On the other hand, $v_y(\,\cdot\,;z)$ solves \eqref{ppl}
subject to the boundary condition $v_y(\,\cdot\,;z)=-1$ on $\partial\mathcal{C}_z$ by Proposition \ref{C1}. Therefore, it admits the standard Feynman-Kac representation (see, e.g., \cite[Ch.\ 5, Sec.\ 7.B]{KS})
\beq
\label{FK}
v_y(x,y;z)=\mathbb{E}\Big[-\,e^{\int_0^{\tau_{\mathcal{C}_z}}\big(\frac{\partial}{\partial y}\mu_2(\tilde{Y}^y_t)-r\big)dt}\Big],
\enq
where $\tau_{\Cc_z}:=\inf\{t\geq 0 \ |\ (X_t^x,\tilde{Y}^y_t)\notin\mathcal{C}_z\},$ and with $\tilde{Y}^y$ solving
\beqs
\begin{cases}
d\tilde{Y}^y_t=\left[\frac{\partial\sigma^2_2}{\partial y}(\tilde{Y}^y_t)+\mu_2(\tilde{Y}^y_t)\right]dt+\sigma_2(\tilde{Y}^y_t)dW_t^2, \qquad t > 0,\\
\tilde{Y}^y_0=y.
\end{cases}
\enqs
Since $r>\frac{\partial\mu_2}{\partial y}$ by Assumption \ref{ass:mu2}, \eqref{FK} implies $v_y(\,\cdot\,;z)> -1$ in $\mathcal{C}_z$, contradicting \eqref{vy<} and concluding the proof.
\ep\\\\
In order to find an upper bound for $y^*(\,\cdot\,;z)$, we now denote
\begin{align}\label{def:F}
F(x,y;z):=c_z(x,z)-\mu_2(y)+ry,\qquad(x,y)\in \overline{Q},
\end{align}
and define
\begin{align}\label{def:vartheta}
\vartheta(x;z):=\inf\{y\in\Ic_2\,|\,F(x,y;z)>0\}\in\overline{\Ic}_2,  \ \ \ x\in\mathcal{I}_1,
\end{align}
with the convention $\inf\emptyset=\overline{y}$. It is worth recalling that \eqref{functional:OS} and standard arguments based on exit times from small subsets of $Q$ give the following inclusion
\begin{align}\label{inclAz}
\mathcal{A}_z\subset L^-_z:=\big\{(x,y)\in Q\,|\,c_z(x,z)\le \mu_2(y)-ry\big\}.
\end{align} 
Then, by Proposition \ref{prop:monot} and by \eqref{inclAz}, we have
\begin{equation}\label{ytheta}
y^*(\,\cdot\,;z)\leq \vartheta(\,\cdot\,;z).
\end{equation}

\begin{Lemma}
\label{rem:nullmeas}
Under Assumption \ref{ass:mu21} and \ref{ass:mu2}, the function $\vartheta(\,\cdot\,;z)$ is nondecreasing and continuous.  
Moreover,  if $\vartheta(x;z)\in\mathcal{I}_2$ then $\vartheta(x;z)$ is the unique solution to the equation $F(x,\cdot;z)=0$ in $\mathcal{I}_2$. Finally one has
\begin{align}\label{ybar}
\big\{(x,y)\in Q\ |\  c_z(x,z)-\mu_2(y)+ry < 0\big\}=\{(x,y)\in Q\ | \ y< \vartheta(x;z)\}.
\end{align}
\end{Lemma}
\textbf{Proof.}
Since $x\mapsto F(x,y;z)$ is nonincreasing (cf.~Assumption \ref{ass:cost}-(ii)) and $y\mapsto F(x,y;z)$ is increasing by Assumption \ref{ass:mu2} and $(x,y)\mapsto F(x,y;z)$ it is not hard to see that $\vartheta(\cdot;z)$ is nondecreasing and right-continuous.

The definition of $\vartheta(\cdot;z)$ and the continuity of $F$ guarantee that if $\vartheta(x;z)\in\mathcal{I}_2$ then $\vartheta(x;z)$ solves $F(x,\cdot;z)=0$ in $\mathcal{I}_2$. Assumption \ref{ass:mu2} then implies that $\vartheta(x;z)$ is actually the unique solution of such equation.

Let us now show that $\vartheta(\,\cdot\,;z)$ is continuous. Take $x_0$ such that $\vartheta(x_0;z)>\underline{y}$ and assume that $\vartheta(x_0-;z)<\vartheta(x_0;z)$. Take a sequence $\{x_n\,,\,n\in\mathbb{N}\}\subset\mathcal{I}_1$ increasing and such that $x_n\uparrow x_0$. One has $F(x_n,\vartheta(x_n;z);z)\ge0$ for all $n\in\mathbb{N}$ and hence in the limit one finds $F(x_0,\vartheta(x_0-;z);z)\ge0\ge F(x_0,\vartheta(x_0;z);z)$ which implies $\vartheta(x_0-;z)\ge\vartheta(x_0;z)$ since $y\mapsto F(x,y;z)$ is increasing.

Clearly \eqref{ybar} follows from the previous properties.
\ep \\

Consider now the class of functions
$$
\mathcal{M}_z:= \{ f:\mathcal{I}_1\rightarrow \overline{\mathcal{I}}_2, \ \mbox{continuous,  nondecreasing and dominated from above by } \vartheta(\,\cdot\,;z)\},
$$
and define
$$\mathcal{D}_f:=\{x\in\mathcal{I}_1 \ | \ f(x)\in\mathcal{I}_2\}, \quad f\in\mathcal{M}_z.$$
Clearly, $\mathcal{M}_z$ is nonempty, as $\vartheta(\,\cdot\,;z) \in \mathcal{M}_z$ by Lemma \ref{rem:nullmeas}.  Moreover $\mathcal{D}_f$ is an open sub-interval (possibly empty) of $\mathcal{I}_1$ due to monotonicity of $f\in\mathcal{M}_z$, that is, 
$$
\mathcal{D}_f =(\underline{x}_f, \overline{x}_f),
$$
where we set
\begin{align}\label{def:xf}
\underline{x}_f:=\inf\{x\in\mathcal{I}_1 \ | \ f(x)>\underline{y}\}, \quad \overline{x}_f:= \sup \{x\in\mathcal{I}_1 \ | \ f(x)<\overline{y}\},
\end{align}
with the conventions $\inf \emptyset = \overline{x}$, $\sup\emptyset=\underline{x}$.
Notice also that by monotonicity of $f \in \mathcal{M}_z$ we have
 $f\equiv \underline{y}$ on $(\underline{x}, \underline{x}_f)$ (if the latter is nonempty) and, analogously, $f\equiv \overline{y}$ on $(\overline{x}_f, \overline{x})$ (if the latter is nonempty).
 Given a function $\hat{y}(\,\cdot\,;z)\in \mathcal{M}_z$, we set
\begin{align}\label{uniq02}
\widehat{H}(x,y;z):=c_z(x,z)\mathds{1}_{\{y>\hat{y}(x;z)\}}-
\big(ry-\mu_2(y)\big)\mathds{1}_{\{y\le\hat{y}(x;z)\}}
\end{align}
and define
\begin{align}\label{defW}
w(x,y;z):=\E\left[\int^\infty_0 e^{-rt}\widehat{H}(X^x_t,Y^y_t;z)dt\right].
\end{align}
Notice that
\begin{align}
\label{eq:subpon02}
\big|w(x,y;z)\big|\le C(z)\big(1+|x|^{\beta}+|y|\big), \qquad\text{for $(x,y)\in Q$,}
\end{align}
by Assumptions \ref{ass:D2}, \ref{ass:cost}, \ref{Ass:r} (cf.~also \eqref{eq:subpon}). Moreover, one can verify that
\beq\label{eq:mg02}
\bigg\{e^{-rt}w(X^x_t,Y^y_t;z)+\int^t_0{e^{-rs}\widehat{H}(X^x_s,Y^y_s;z)ds},\:t\ge0\bigg\}\:\:\:\text{is an $(\Fc_t)$-martingale}
\enq
and the family $\big\{e^{-r\tau}w(X^x_\tau,Y^y_\tau;z)\,,\,\tau\in\Tc\big\}$ is uniformly integrable.

To simplify notations, from now on we set
\begin{equation}
\label{notationhat}
\hat{x}:= \overline{x}_{\hat{y}(\cdot;z)}, \qquad \check{x}:=\underline{x}_{\hat{y}(\cdot;z)}, \qquad \hat{\mathcal{D}}_z:= {\mathcal{D}}_{\hat{y}(\cdot;z)},
\end{equation}
and
\begin{equation}
\label{notationstar}
{x}^*:= \overline{x}_{{y}^*(\cdot;z)}, \qquad {x}_*:=\underline{x}_{{y}^*(\cdot;z)}, \qquad {\mathcal{D}}^*_z:= {\mathcal{D}}_{{y}^*(\cdot;z)}.
\end{equation}

We can now state the main result of this section. We use arguments inspired by \cite[Sec.~25]{PeskShir} and references therein.
\begin{Theorem}
\label{teo:uniqueness}
Let Assumptions \ref{ass:mu21}, \ref{ass:densitiesderivative} and \ref{ass:mu2} hold. Assume that $\mathcal{C}_z\neq \emptyset$ and $\mathcal{A}_z\neq \emptyset$. Then $y^*(\,\cdot\,;z)$ is the unique function $y(\,\cdot\,;z)\in \mathcal{M}_z$ with $\mathcal{D}_{y(\cdot;z)}\neq \emptyset$ and such that for each $x\in \mathcal{D}_{y(\cdot;z)}$ it holds
\begin{align}
\label{eq:intan2}
-y(x;z) = &\int_0^\infty e^{-rt} \bigg[\int^{\overline{x}}_{\underline{x}}p_1(t,x,\xi)c_z(\xi,z)\bigg(\int^{\overline{y}}_{y(\xi;z)} \ p_2(t,y(x;z),\eta) d\eta\bigg)d\xi\bigg]dt\\
&  -\,\int_0^\infty e^{-rt} \bigg[\int^{\overline{x}}_{\underline{x}} p_1(t,x,\xi)\bigg(\int^{y(\xi;z)}_{\underline{y}} (r\eta-\mu_2(\eta )) p_2(t,y(x;z),\eta)d\eta\bigg)d\xi\bigg]dt.\nonumber
\end{align}
\end{Theorem}
\textbf{Proof.}
\emph{Existence.}
First of all, we observe that $y^*(\,\cdot\,;z)\in\mathcal{M}_z$ by Propositions \ref{prop:y*}, \ref{continuousboundary}, and \eqref{ytheta}.
The fact that $y^*(\,\cdot\,;z)$ solves \eqref{eq:intan2} for each  $x\in \mathcal{D}^*_z$
follows by evaluating both sides of \eqref{v:rep1} at points of the boundary $(x,y^*(x;z)) \in \partial\mathcal{A}_z$, which yields
\begin{eqnarray}
\label{eq:int}
-y^*(x;z)& \hspace{-0.25cm} = \hspace{-0.25cm} & \int_0^\infty e^{-rt}\mathbb{E}\Big[c_z(X^x_t,z)\mathds{1}_{\{Y^{y^*(x;z)}_t>y^*(X^x_t;z)\}}\Big]dt \\
& & - \,\int_0^\infty e^{-rt}\mathbb{E}\Big[(rY^{y^*(x;z)}_t-\mu_2(Y^{y^*(x;z)}_t))\mathds{1}_{\{Y^{y^*(x;z)}_t\leq y^*(X^x_t;z)\}}\Big]dt.\nonumber
\end{eqnarray}
From \eqref{eq:int} and by Assumption \ref{ass:laws},  we see that $y^*(\,\cdot\,;z)$ solves \eqref{eq:intan2}.

\medskip
\emph{Uniqueness.}
Recall \eqref{notationhat} and \eqref{notationstar}.
Let $\hat{y}(\,\cdot\,;z)\in\mathcal{M}_z$ be such that $\hat{\mathcal{D}_z}\neq \emptyset$ and  solving  \eqref{eq:intan2} on $\hat{\mathcal{D}_z}$.
 We need to show that $\hat{y}(\,\cdot\,;z)\equiv y^*(\,\cdot\,;z)$.

\smallskip
\emph{Step 1.} Here we show that  $\hat{y}(\,\cdot\,;z)\ge y^*(\,\cdot\,;z)$. We distinguish two cases: when $\Dc^*_z\cap\hat{\Dc}_z\neq\emptyset$ and when $\Dc^*_z\cap\hat{\Dc}_z = \emptyset$. Notice that in general $\Dc^*_z\cap\hat{\Dc}_z=(x_*\vee \check{x},x^*\wedge\hat{x})$.

\emph{Case $\Dc^*_z\cap\hat{\Dc}_z\neq\emptyset$.} First, we show that $\hat{y}(\,\cdot\,;z)\geq y^*(\,\cdot\,;z)$ on $\Dc^*_z\cap\hat{\Dc}_z$ and later we will prove it on $\mathcal{I}_1\setminus (\Dc^*_z\cap\hat{\Dc}_z)$.
 Assume, by contradiction, that $\hat{y}(x;z)< y^*(x;z)$ for some $x\in\Dc^*_z\cap\hat{\Dc}_z$, take $y<\hat{y}(x;z)$ and set $\sigma=\sigma(x,y,z):=\inf\big\{t\ge0\,|\,Y^y_t\ge y^*(X^x_t;z)\big\}$.
Then, from \eqref{eq:mg} and \eqref{eq:mg02}, it follows (up to usual localization arguments as in 
Lemma \ref{lemma:w}) that
\begin{align}
\vspace{+10pt}
&\E\left[e^{-r\sigma}v(X^x_\sigma,Y^y_\sigma;z)\right]=v(x,y;z)+\E\left[\int^\sigma_0{e^{-rt}
\big(rY^y_t-\mu_2(Y^y_t)\big)\,dt}\right],\label{uniq04-a}\\
&\E\left[e^{-r\sigma}w(X^x_\sigma,Y^y_\sigma;z)\right]=w(x,y;z)-\E\left[\int^\sigma_0{e^{-rt}
\hat{H}(X^x_t,Y^y_t;z)\,dt}\right].\label{uniq04-b}
\end{align}
Lemma \ref{lemma:w} in Appendix \ref{app:os1d} ensures that $v\ge w$ everywhere and that $w(x,y;z)=v(x,y;z)=-y$, since $y<\hat{y}(x;c)<y^*(x;z)$ (cf.~\eqref{ww}). Then, subtracting \eqref{uniq04-b} from \eqref{uniq04-a}, one has
\begin{align}\label{uniq05}
0\le& \ \E\left[\int^\sigma_0{e^{-rt}\left[
\big(rY^y_t-\mu_2(Y^y_t)\big)+\hat{H}(X^x_t,Y^y_t;z)\right]\,dt}\right]\nonumber\\
=&\ \E\left[\int^\sigma_0{e^{-rt}
\left[c_z(X^x_t,z)-\big(\mu_2(Y^y_t)-rY^y_t\big)\right]\mathds{1}_{\{\hat{y}(X^x_t;z)<Y^y_t<y^*(X^x_t;z)\}}\,dt}\right].
\end{align}
Notice that the continuity of trajectories of $(X^x,Y^y)$ and the continuity of $y^*(\,\cdot\,;z)$ give $\sigma>0$ $\P$-a.s. Moreover, from the continuity of $y^*(\,\cdot\,;z)$ and $\hat{y}(\,\cdot\,;z)$ one gets that the set $\big\{(x,y)\in Q\ |\ \hat{y}(x;z)<y<y^*(x;z)\big\}$ is open and not empty. These facts, combined with the fact that $y^*(\,\cdot\,;z)\leq \vartheta(\,\cdot\,;z)$ and with \eqref{ybar}, imply that the last expression in \eqref{uniq05} must be strictly negative and we reach a contradiction. 
Therefore 
\begin{equation}\label{aa1}
\hat{y}(x;z)\ge y^*(x;z),\quad \text{for all $x\in \Dc^*_z\cap\hat{\Dc}_z=(\check{x}\vee x_*\,,\,\hat{x}\wedge x^*)$.}
\end{equation}

Now we show that  $\hat{y}(\,\cdot\,;z)\ge y^*(\,\cdot\,;z)$ on $\mathcal{I}_1\setminus(\Dc^*_z\cap\hat{\Dc}_z)$, if $(\Dc^*_z\cap\hat{\Dc}_z)\neq \emptyset$.
By \eqref{aa1} and continuity of $\hat{y}(\,\cdot\,;z)$ and $y^*(\,\cdot\,;z)$, we deduce that the inequality \eqref{aa1} also holds at the endpoints of the interval, i.e.
\begin{equation}
\label{a0}
\hat{y}(x_*\vee \check{x};z)\geq {y}^*(x_*\vee \check{x};z) \quad\text{and}\quad \hat{y}(x^*\wedge \hat{x};z)\geq y^*(x^*\wedge \hat{x};z).
\end{equation}  
It follows from the definition of $\check{x}$ and continuity of $\hat{y}(\,\cdot\,;z)$ that $\underline{y}=\hat{y}(\check{x};z)$. So, if we argue by contradiction and assume $\check{x}>x_*$, then by definition of $x_*$ we get $y^*(\check{x};z)>\underline{y}$. The latter and the first inequality in \eqref{a0} imply $\underline{y}=\hat{y}(\check{x};z)\ge y^*(\check{x};z)>\underline{y}$, hence a contradiction. Thus, we conclude $\check{x}\le x_*$. 

By an analogous argument applied to the second inequality in \eqref{a0} we also obtain $x^*\ge \hat{x}$ and therefore 
\begin{align*}
\Dc^*_z\cap\hat{\Dc}_z=(\check{x}\vee x_*\,,\,\hat{x}\wedge x^*)=(x_*\,,\,\hat{x}).
\end{align*}
By monotonicity and continuity of $\hat{y}(\,\cdot\,;z)$ and $y^*(\,\cdot\,;z)$, and by definition of $\hat{x}$ and $x_*$ we have
\begin{align*}
\text{$y^*(x;z)=\underline{y}$ for $x\le x_*$ and $\hat{y}(x;z)=\overline{y}$ for $x\ge \hat{x}$.}
\end{align*}
On the other hand $y^*(x;z)\le \overline{y}$ and $\hat{y}(x;z)\ge \underline{y}$ for all $x\in\mathcal{I}_1$ and therefore $\hat{y}(\,\cdot\,;z)\ge y^*(\,\cdot\,;z)$ on $\mathcal{I}_1\setminus(\Dc^*_z\cap\hat{\Dc}_z)$ as claimed.

\smallskip
\emph{Case  $\Dc^*_z\cap\hat{\Dc}_z=\emptyset$}. By monotonicity of $y^*(\,\cdot\,;z)$ and $\hat{y}(\,\cdot\,;z)$, one has either $\hat{x}\le x_*$ or $\check{x}\ge x^*$. If $\hat{x}\le x_*$, then $\hat{y}(\,\cdot\,;z)\ge y^*(\,\cdot\,;z)$ on $\Ic_1$; if $\check{x}\ge x^*$, we can use the same arguments as above to find $\check{x}=\overline{x}$, which contradicts the assumption that   $\hat{\mathcal{D}}_z\neq \emptyset$.

\medskip
\emph{Step 2.} Here we show that  $\hat{y}(\,\cdot\,;z)\le y^*(\,\cdot\,;z)$.
Assume, by contradiction,  that there exists $x\in\Ic_1$ such that $\hat{y}(x;z)> y^*(x;z)$. Take $y\in(y^*(x;z)\,,\,\hat{y}(x;z))$ and consider the stopping time $\tau^* = \tau^*(x,y;z):= \inf\{t\geq 0\ | \ Y^y_t\leq y^*(X^x_t;z)\}$. This is the first optimal stopping time for the problem \eqref{P:OS}, as it is the first entry time in the stopping region $\mathcal{A}_z$ (cf.\ \eqref{superarmonic2} and \eqref{ACz}). As in \emph{Step} 1 above, \eqref{eq:mg} and \eqref{eq:mg02} give
\begin{align}
\vspace{+10pt}
&\E\left[e^{-r{\tau^*}}v(X^x_{\tau^*},Y^y_{\tau^*};z)\right]=v(x,y;z)-\E\left[\int^{\tau^*}_0{e^{-rt}
c_z(X^x_t,z)\,dt}\right],\label{uniq06-a}\\
&\E\left[e^{-r{\tau^*}}w(X^x_{\tau^*},Y^y_{\tau^*};z)\right]=w(x,y;z)-\E\left[\int^{\tau^*}_0{e^{-rt}
\widehat{H}(X^x_t,Y^y_t;z)\,dt}\right].\label{uniq06-b}
\end{align}
 By using \eqref{uiboundv} and standard localization argument, 
we obtain  $\E\big[e^{-r{\tau^*}}v(X^x_{\tau^*},Y^y_{\tau^*};z)\big]=-\E\big[e^{-r{\tau^*}}Y^y_{\tau^*}\big]$. On the other hand, we know from \emph{Step} 1 above that $\hat{y}(\,\cdot\,;z)\ge y^*(\,\cdot\,;z)$, hence $\E\big[e^{-r{\tau^*}}w(X^x_{\tau^*},Y^y_{\tau^*};z)\big]=-\E\big[e^{-r{\tau^*}}Y^y_{\tau^*}\big]$ by \eqref{eq:ww00}, \eqref{ww}, the fact that $\overline{y}$ is a natural boundary point and by localization arguments as in the proof of Lemma \ref{lemma:w}. Taking also into account that $v\ge w$ (cf.~Lemma \ref{lemma:w}) and subtracting \eqref{uniq06-b} from \eqref{uniq06-a} we obtain
\beq\label{uniq07}
0 \hspace{-0.25cm} & \ge & \hspace{-0.25cm} \E\left[\int^{\tau^*}_0{e^{-rt}\left(
\widehat{H}(X^x_t,Y^y_t;z)-c_z(X^x_t,z)\right)\,dt}\right]\nonumber\\
\hspace{-0.25cm} & = & \hspace{-0.25cm} -\E\left[\int^{\tau^*}_0{e^{-rt}\left(
c_z(X^x_t,z)+(rY^y_t-\mu_2(Y^y_t))\right)\mathds{1}_{\{y^*(X^x_t;z)<Y^y_t<\hat{y}(X^x_t;z)\}}\,dt}\right].
\enq
Now $\tau^* > 0$ $\P$-a.s.\ by continuity of trajectories of $(X^x,Y^y)$ and of $y^*(\,\cdot\,;z)$. Moreover, the set $\big\{(x,y)\in Q\,|\,y^*(x;z)<y<\hat{y}(x;z)\big\}$ is open in $Q$ and not empty, by continuity of $y^*(\,\cdot\,;z)$ and $\hat{y}(\,\cdot\,;z)$. Since by assumption $\hat{y}(\,\cdot\,;z)\le \vartheta(\,\cdot\,;z)$, these facts together with \eqref{ybar} imply that the last term in \eqref{uniq07} must be strictly positive thus leading to a contradiction. Hence, $\hat{y}(\,\cdot\,;z)\le y^*(\,\cdot\,;z)$.
\hfill$\square$

\begin{Remark}
\label{rem-inteq}
It is interesting to formulate \eqref{eq:intan2} in the canonical Fredholm form as there exists a wide literature on numerical methods for this kind of nonlinear integral equations. One can rewrite \eqref{eq:intan2} as
\begin{align}
\label{eq:int-rem}
-y(x;z) = &\int_0^\infty e^{-rt} \bigg[\int^{\overline{x}}_{\underline{x}}p_1(t,x,\xi)c_z(\xi,z)d\xi\bigg]dt\\
&  -\,\int_0^\infty e^{-rt} \bigg[\int^{\overline{x}}_{\underline{x}} p_1(t,x,\xi)\bigg(\int^{y(\xi;z)}_{\underline{y}} (c_z(\xi,z) + r\eta-\mu_2(\eta )) p_2(t,y(x;z),\eta)d\eta\bigg)d\xi\bigg]dt.\nonumber
\end{align}
Then, defining
$$K(x,\xi,\alpha,\beta,z):=\int_0^\infty e^{-rt} p_1(t,x,\xi) \Big( \int_{\underline{y}}^{\beta} (c_z(\xi,z) + r\eta-\mu_2(\eta )) p_2(t,\alpha,\eta)d\eta\Big) dt,$$
$$f(x;z):=\int_0^\infty e^{-rt} \bigg[\int^{\overline{x}}_{\underline{x}}p_1(t,x,\xi)c_z(\xi,z)d\xi\bigg]dt,$$
and, applying Fubini's Theorem, 
one finds that \eqref{eq:intan2} takes the form
\begin{equation}
\label{eq:int2}
-y(x;z) = f(x;z) - \int_{\underline{x}}^{\overline{x}} K(x,\xi,y(x;z),y(\xi;z),z) d\xi.
\end{equation}
The latter is a nonlinear Fredholm integral equation of second kind, possibly singular if $\Ic_1=(\underline{x}, \overline{x})$ is unbounded (see, e.g., \cite{Delves} or \cite{Hackbusch}). 
A survey of numerical methods for equations of this kind may be found in \cite{Atkinson} (see also classical textbooks like \cite{Baker} and \cite{Delves}). These methods can be used to solve our equation \eqref{eq:int2}. However, since they are certainly non trivial, we believe that such numerical computation falls outside the scopes of our work.
\end{Remark}

Regarding the assumptions $\mathcal{C}_z\neq \emptyset$ and $\mathcal{A}_z\neq \emptyset$ in Theorem \ref{teo:uniqueness}, we provide the following characterization.
\begin{Proposition}\label{prop:notempty}
\begin{itemize}
\item[1.] The continuation set $\Cc_z$ is not empty if and only if the set
\begin{align}\label{def:Lplus}
L^+_z:=\big\{(x,y)\in Q\,|\,c_z(x,z)-\mu_2(y)+ry>0\big\}
\end{align}
is not empty.
\item[2.] The stopping set $\Ac_z$ is not empty if and only if 
\begin{align}\label{eq:limnoempty}
\lim_{x\uparrow \overline{x}}\E\left[\int_0^\infty{e^{-rt}c_z(X^x_t,z)dt}\right]<-\underline{y}.
\end{align}
\end{itemize}
\end{Proposition}
\textbf{Proof.}
For the first claim notice that $L^+_z\subset \Cc_z$ (cf.~also \eqref{inclAz}) so that $L^+_z\neq\emptyset\Rightarrow\Cc_z\neq\emptyset$. To prove the reverse implication it suffices to observe that, by using \eqref{integroparti} into \eqref{functional:OS}, if $L^+_z=\emptyset$ then any stopping rule would produce a payoff smaller or equal than the one of immediate stopping and therefore $\Cc_z=\emptyset$.

For the second claim we observe that
$$\Ac_z=\emptyset\ \Longleftrightarrow \ \Cc_z=Q\ \Longleftrightarrow\ \tau^*=+\infty \ \P-a.s.\ \forall (x,y)\in Q \ \Longleftrightarrow\  v(x,y;z)>-y\ \  \forall (x,y)\in Q.$$ Hence, $\mathcal{A}_z= \emptyset$ if and only if
\begin{align}\label{eq:vnostop}
v(x,y;z)=\E\left[\int_0^\infty{e^{-rt}c_z(X^x_t,z)dt}\right]>-y\quad\forall (x,y)\in Q.
\end{align}
Then, \eqref{eq:limnoempty} implies that $\Ac_z\neq\emptyset$. Conversely, if $\Ac_z\neq \emptyset$, then there exists a point $(x,y)\in Q$ such that stopping at once is more profitable than (for instance) never stopping. For such a point
\begin{align}\label{eq:vnostop2}
0=y+v(x,y;z)\ge y+\E\left[\int_0^\infty{e^{-rt}c_z(X^x_t,z)dt}\right].
\end{align}
Since $y> \underline{y}$ and $c_z(\,\cdot\,,z)$ is nonincreasing (cf.~Assumption \ref{ass:cost}-(ii)), 
then \eqref{eq:limnoempty} must hold.
\ep

\medskip

In principle Theorem \ref{teo:uniqueness} fully characterizes the optimal boundary of problem \eqref{P:OS}, but it has the drawback that the region $\mathcal{D}^*_z = ({x}_*, {x}^*)$, with ${x}_*$ and ${x}^*$ as in \eqref{notationstar},
 is defined implicitly. For the purpose of numerical evaluation of \eqref{eq:intan2} it would be helpful to know $\Dc^*_z$, in advance rather than computing it at the same time as $y^*(\,\cdot\,;z)$.
Recall \eqref{def:xf} and define
\begin{align}
{\theta}_*:=\underline{x}_{\vartheta(\cdot;z)}=\inf\big\{x\in\Ic_1\,|\,\vartheta(x;z)>\underline{y}\big\},\quad{\theta}^*:=\overline{x}_{\vartheta(\cdot;z)}=\sup\big\{x\in\Ic_1\ |\ \vartheta(x;z)<\overline{y}\big\},\label{def:thetaz}
\end{align}
with the convention $\inf\emptyset=\overline{x}$,  $\sup\emptyset=\underline{x}$. Since $y^*(\,\cdot\,;z)\leq \vartheta(\,\cdot\,;z)$, we have $x_*\geq \theta_*$ and $x^*\geq \theta^*$.
To characterize $x_*$ we will make use of the following algebraic equation
\begin{align}\label{eq:x*1}
-\underline{y}=\int_0^\infty{e^{-rt}\Big(\int^x_{\underline{x}}{p_1(t,x;\xi)c_z(\xi,z)d\xi}-r\underline{y}\int_x^{\overline{x}}p_1(t,x,\xi)d\xi\Big)dt}.
\end{align}
Similarly, if $\overline{y}<+\infty$, a characterization of $x^*$ will be given in terms of the algebraic equation
\begin{align}\label{eq:x*2}
-\overline{y}=\int_0^\infty{e^{-rt}\Big(\int^x_{\underline{x}}{p_1(t,x;\xi)c_z(\xi,z)d\xi}-r\overline{y}\int_x^{\overline{x}}p_1(t,x,\xi)d\xi\Big)dt}.
\end{align}

\begin{Proposition}\label{prop:shapeD}
Let Assumptions \ref{ass:mu21}, \ref{ass:densitiesderivative}, \ref{ass:mu2} hold.  Let $\mathcal{C}_z\neq \emptyset$ and $\mathcal{A}_z\neq \emptyset$.
Then
\begin{itemize}
\item[1.] $x_*\in \Ic_1$ if and only if \eqref{eq:x*1} has a unique solution $\tilde{x}\in(\theta_*,\overline{x})$, and, in this case, $x_*=\tilde{x}$. Otherwise, we have $x_*=\underline{x}$.
\item[2.] If $\overline{y}<+\infty$, then $x^*\in \Ic_1$ if and only if \eqref{eq:x*2} has a unique solution $\tilde{x}'\in(\theta^*,\overline{x})$, and, in this case, $x^*=\tilde{x}'$. Otherwise, we have   $x^*=\overline{x}$.
\item[3.] If $\overline{y}=+\infty$ and there exists $\lambda>0$ such that $r-\frac{\partial \mu_2}{\partial y}\ge \lambda$ on $\Ic_2$, then $x^*=\overline{x}$.
\end{itemize}
\end{Proposition}
\textbf{Proof.}
1. Existence and uniqueness of a solution of \eqref{eq:x*1} $(\theta_*,\overline{x})$ is discussed in Appendix \ref{A2}.

\emph{Proof of} $\Rightarrow$.
Take a sequence $\{x_n,\,n\in\mathbb{N}\}\subset\Ic_1$ such that $x_n\downarrow x_*$ and notice that by Theorem \ref{teo:uniqueness} we have for every $n\in\mathbb{N}$
\begin{align}\label{eq:intan3}
-y^*(x_n;z) = &\int_0^\infty e^{-rt} \bigg[\int^{\overline{x}}_{\underline{x}}p_1(t,x_n,\xi)c_z(\xi,z)\bigg(\int^{\overline{y}}_{y^*(\xi;z)} \ p_2(t,y^*(x_n;z),\eta) d\eta\bigg)d\xi\bigg]dt\\
&  -\,\int_0^\infty e^{-rt} \bigg[\int^{\overline{x}}_{\underline{x}} p_1(t,x_n,\xi)\bigg(\int^{y^*(\xi;z)}_{\underline{y}} (r\eta-\mu_2(\eta )) p_2(t,y^*(x_n;z),\eta)d\eta\bigg)d\xi\bigg]dt.\nonumber
\end{align}
We aim to take limits of \eqref{eq:intan3} as $n\uparrow \infty$.
For the left hand-side of \eqref{eq:intan3} we have $y^*(x_n;z)\downarrow \underline{y}$, by continuity of $y^*(\,\cdot\,;z)$ and definition of $x_*$.
On the other hand, taking into account that $y^*(\,\cdot\,;z)=\underline{y}$ for $\xi\leq x_*$, the first term of the right-hand side of \eqref{eq:intan3} can be written as
\begin{align}
\label{sse}
&\int_0^\infty e^{-rt} \bigg[\int^{\overline{x}}_{\underline{x}}p_1(t,x_n,\xi)c_z(\xi,z)\bigg(\int^{\overline{y}}_{y^*(\xi;z)} \ p_2(t,y^*(x_n;z),\eta) d\eta\bigg)d\xi\bigg] dt\\
=&
\int_0^\infty e^{-rt} \bigg[\int^{x_*}_{\underline{x}}p_1(t,x_n,\xi)c_z(\xi,z)d\xi \nonumber \\
&\hspace{2cm}+ \int^{\overline{x}}_{x_*}p_1(t,x_n,\xi)c_z(\xi,z)\bigg(\int^{\overline{y}}_{\underline{y}} \ \mathds{1}_{\{\eta > y^*(\xi;z)\}} p_2(t,y^*(x_n;z),\eta) d\eta\bigg)d\xi\bigg]dt.\nonumber
\end{align}
\noindent Now notice that:

\emph{(i)}\, for any $t>0$ the sequence of probability measures with densities $\{p_1(t,x_n,\xi),\,n \in \mathbb{N}\}$ on $\mathcal{I}_1$ converges pointwisely to $p_1(t,x_*,\xi)d\xi $ by Assumption \ref{ass:laws};

\emph{(ii)}\,as $\underline{y}$ is non-entrance (since natural), for any given and fixed $t>0$ and $z \in \mathbb{R}^+$ the sequence of probability measures with densities $\{p_2(t,y^*(x_n;z),\eta),\,n \in \mathbb{N}\}$ on $\mathcal{I}_2$ converges weakly to the Dirac's delta measure $\delta_{\underline{y}}(\eta)$ (see Section \ref{A1});

\emph{(iii)}\, for every $\xi>x_*$, the function $\mathcal{I}_2\rightarrow \mathbb{R}$,
$\eta\mapsto c_z(\xi,z)\mathds{1}_{\{\eta > y^*(\xi;z)\}} \equiv 0$ $\delta_{\underline{y}}$-a.e.

Then, taking into account \emph{(i)-(iii)} we can apply  Portmanteau Theorem to the integral with respect to $d\eta$ in the right hand side of \eqref{sse} and dominated convergence to the one with respect to $d\xi$ to obtain
$$
\lim_{n\rightarrow+\infty} \int^{\overline{x}}_{\underline{x}}p_1(t,x_n,\xi)c_z(\xi,z)\bigg(\int^{\overline{y}}_{y^*(\xi;z)} \ p_2(t,y^*(x_n;z),\eta) d\eta\bigg)d\xi= \int^{x_*}_{\underline{x}}p_1(t,x_*,\xi)c_z(\xi,z)d\xi
$$
Finally, a further application of dominated convergence to the integral with respect to $dt$, gives
\begin{align*}
\lim_{n\rightarrow+\infty}&\int_0^{\infty} e^{-rt}\left[ \int^{\overline{x}}_{\underline{x}}p_1(t,x_n,\xi)c_z(\xi,z)\bigg(\int^{\overline{y}}_{y^*(\xi;z)} \ p_2(t,y^*(x_n;z),\eta) d\eta\bigg)d\xi\right] dt
\\
=&\int_0^\infty e^{-r t} \left[ \int^{x_*}_{\underline{x}}p_1(t,x_*,\xi)c_z(\xi,z)d\xi\right] dt.
\end{align*}
Similar arguments can be applied to the second term of the right-hand side of \eqref{eq:intan3}. In fact for $\xi>x_*$ the map $\eta\mapsto(r\eta-\mu_2(\eta))\mathds{1}_{\{\eta\le y^*(\xi;z)\}}$ is bounded on $\overline{\Ic}_2$ and it is continuous at $\underline{y}$. Moreover $(r\eta-\mu_2(\eta))\mathds{1}_{\{\eta\le y^*(\xi;z)\}}=r\underline{y}-\mu_2(\underline{y})$, $\delta_{\underline{y}}$-a.e.

\emph{Proof of} $\Leftarrow$. Assume now that $\theta_*<\overline{x}$ and that $\tilde{x}\in(\theta_*,\overline{x})$ uniquely solves \eqref{eq:x*1}. It is proven in Appendix \ref{app:os1d}, Section \ref{A2}, that $\tilde{x}$ is the optimal boundary of the one-dimensional optimal stopping problem
\begin{align}
\label{eq:undv}
\underline{v}(x;z):=&\sup_{\tau\in\Tc}\E\left[\int_0^\tau{e^{-rt}c_z(X^x_t,z)dt}-\underline{y}e^{-r\tau}\right],
\end{align}
 hence that $\underline{\Ac}_z:=\{x\in\Ic_1\,|\,\underline{v}(x;z)=-\underline{y}\}=\{x\in \Ic_1\,|\,x\ge\tilde{x}\}$. By arguments as in the proof of Proposition \ref{vcont} we have $\underline{v}(x;z)=\lim_{y\downarrow\underline{y}}v(x,y;z)$. Moreover, $0<\underline{v}(x;z)+\underline{y}\le v(x,y;z)+y$ for all $(x,y)\in(\underline{x},\tilde{x})\times\Ic_2$, by monotonicity of $y\mapsto v(x,y;z)+y$ (cf.\ Proposition \ref{prop:monot}); hence $x_*\ge\tilde{x}>\underline{x}$. Also $x_*<\overline{x}$, since otherwise $\Ac_z=\emptyset$, thus contradicting the assumption $\Ac_z\neq \emptyset$.
Therefore, $x_*\in\Ic_1$ and, by the arguments of the first part of this proof, $x_*$ solves \eqref{eq:x*1}. Since such solution is unique, it must be $\tilde{x}=x_*$.

\medskip

2. The proof of this second claim works thanks to arguments similar to the ones employed for the first one. One has to consider, in place of \eqref{eq:undv}, the optimal stopping problem
\begin{align*}
\overline{v}(x;z):=&\sup_{\tau\in\Tc}\E\left[\int_0^\tau{e^{-rt}c_z(X^x_t,z)dt}-\overline{y}e^{-r\tau}\right].
\end{align*}

3. The further assumption guarantees that $\vartheta(\,\cdot\,;z)<+\infty$ on $\Ic_1$ and the claim follows.
\ep
\medskip

\begin{Remark}
\label{boundaries1d}
Despite their rather involved definition, $x_*$ and $x^*$ have a quite clear probabilistic interpretation. In fact, they are the free-boundaries of the optimal stopping problems
\begin{align*}
\underline{v}(x;z):=&\sup_{\tau\in\Tc}\E\left[\int_0^\tau{e^{-rt}c_z(X^x_t,z)dt}-\underline{y}e^{-r\tau}\right],\quad
\overline{v}(x;z):=\sup_{\tau\in\Tc}\E\left[\int_0^\tau{e^{-rt}c_z(X^x_t,z)dt}-\overline{y}e^{-r\tau}\right],
\end{align*}
respectively, with $\underline{v}(\,\cdot\,;z)=\lim_{y\downarrow \underline{y}}v(\,\cdot\,,y;z)$ and $\overline{v}(\,\cdot\,;z)=\lim_{y\uparrow \overline{y}}v(\,\cdot\,,y;z)$.
\end{Remark}

\section{The Optimal Control}
\label{sec:opt}

In this section we characterize the optimal control $\nu^*$ of \eqref{OCP} by showing that it is optimal to exert the minimal effort needed to reflect the (optimally controlled) state process $Z^{z,\nu^*}$ at a (random) boundary intimately connected to $y^*$ of Theorem \ref{teo:uniqueness}. 

\subsection{The action/inaction regions}
Define
\begin{align}
\label{regionsCA0}
\mathcal{C}:=\{(x,y,z)\in\mathcal{O} \ | \ v(x,y;z)>-y\}\quad\text{and}\quad \mathcal{A}:=\{(x,y,z)\in\mathcal{O} \ | \ v(x,y;z)=-y\}.
\end{align}
The sets $\Cc$ and $\Ac$ are respectively the candidate inaction region and the candidate action region for the control problem \eqref{OCP}.

\begin{Remark}
We notice that from the connection proved in \cite{BK} we expect $V_z=v$ and
\begin{align}
\label{regionsCA}
\mathcal{C}=\{(x,y,z)\in\mathcal{O} \ | \ V_z(x,y,z)>-y\}, \ \ \ \mathcal{A}=\{(x,y,z)\in\mathcal{O} \ | \ V_z(x,y,z)=-y\}.
\end{align}
Intuitively, $\mathcal{A}$ is the region in which it is optimal to invest immediately, whereas $\mathcal{C}$ is the region in which it is profitable to delay the investment option.
\end{Remark}

Throughout this section all the assumptions made so far will be standing assumptions, i.e.\ Assumptions \ref{ass:D2}, \ref{ass:cost}, \ref{Ass:r}, \ref{ass:ui}, \ref{ass:laws}, \ref{ass:mu21}, \ref{ass:densitiesderivative} and \ref{ass:mu2} hold and we will not repeat them in the statement of the next results.

It immediately follows from the fact that $c_z(x,\cdot)$ is nondecreasing for each $x\in \mathcal{I}_1$ that
\begin{Proposition}
\label{vdecrinz}
The function $z\mapsto v(x,y;z)$ is nondecreasing for every $(x,y)\in Q$.
\end{Proposition}
The nondecreasing property of $z \mapsto v(x,y;z)$ implies that for fixed $(x,y)\in Q$ the region $\mathcal{A}$ is below $\mathcal{C}$, and we define the boundary between these two regions by
\beq
\label{defzstar}
z^*(x,y):=\inf\{z\in\mathbb{R}^+\ | \ v(x,y;z)>-y\},
\enq
with the convention $\inf \emptyset =\infty$. Then \eqref{regionsCA0} can be equivalently written as
\begin{align}
\label{regionsCA1}
\mathcal{C}=\{(x,y,z)\in\mathcal{O} \ | \ z > z^*(x,y)\},\ \ \
\mathcal{A}=\{(x,y,z)\in\mathcal{O} \ | \ z \leq z^*(x,y)\}.
\end{align}
We can also easily observe from \eqref{ystar} and \eqref{defzstar} and from the nondecreasing property of $z \mapsto v(x,y;z)$ and of $y \mapsto v(x,y;z) + y$ (cf.\ Proposition \ref{vdecrinz} and Proposition \ref{prop:monot}, respectively) that
\begin{eqnarray}
\label{equivalences}
z \ > \ z^*(x,y)\ \Longleftrightarrow  \ v(x,y;z) \ > \ - y \ \Longleftrightarrow \ y \ > \ y^*(x;z), \quad (x,y,z) \in \mathcal{O}.
\end{eqnarray}
Hence, for any $x\in \mathcal{I}_1$, $z^*$ of \eqref{defzstar} can be seen as the pseudo-inverse of the nonincreasing (cf.\ Proposition \ref{prop:y*}) function $z \mapsto y^*(x;z)$; that is,
\beq
\label{z*y*}
z^*(x,y)=\inf\{ z\in \mathbb{R}^+ \ |\ y> y^*(x;z)\}, \quad (x,y) \in Q.
\enq
It thus follows that the characterization of $y^*$ of Theorem \ref{teo:uniqueness} is actually equivalent to a complete characterization of $z^*$ thanks to \eqref{z*y*}.

Set $$\overline{z}(x,y) := \inf\{z\in\mathbb{R}^+\ | \ c_z(x,z)-\mu_2(y)+ry>0\}, \ \ \ (x,y)\in Q,$$
with the usual convention $\inf\emptyset =\infty$, and recall $\vartheta(x;z)$ of Lemma \ref{rem:nullmeas}. Then the nondecreasing property of
$z \mapsto c_z(x,z)-\mu_2(y)+ry$ and of $y \mapsto c_z(x,z)-\mu_2(y)+ry$ (cf.\ Assumption \ref{ass:cost} and Assumption \ref{ass:mu21}, respectively) implies that
$$z \ > \ \overline{z}(x,y)\ \Longleftrightarrow  \ c_z(x,z)-\mu_2(y)+ry>0 \ \Longleftrightarrow \ y \ > \ \vartheta(x;z), \quad (x,y,z) \in \mathcal{O},$$
and therefore that
\beq
\label{zetabarybar}
\overline{z}(x,y) = \inf\{z\in\mathbb{R}^+\ | \ y > \vartheta(x;z) \}.
\enq
\begin{Proposition}
\label{prop:z*}
One has
\begin{enumerate}
\item $z^*\leq \overline{z}$ over $Q$.
\item $z^*(\,\cdot\,,y)$ is nondecreasing for each $y\in\mathcal{I}_2$ and $z^*(x,\,\cdot\,)$ is nonincreasing for each $x\in\mathcal{I}_1$.
\item $z^*(\,\cdot\,,y)$ is right-continuous for each $y\in\mathcal{I}_2$ and $z^*(x,\,\cdot\,)$ is left-continuous for each $x\in\mathcal{I}_1$.
\item $(x,y) \mapsto z^*(x,y)$ is upper-semicontinuous.
\end{enumerate}
\end{Proposition}
\textbf{Proof.}
\emph{1.} It follows by \eqref{z*y*}, \eqref{zetabarybar} and \eqref{inclAz}.

\emph{2.} The first claim follows from the fact that $v(\,\cdot\,,y;z)$ is nonincreasing for each $y\in \mathcal{I}_2$, $z\in \mathbb{R}^+$, by Proposition \ref{prop:v}; the fact that $y\mapsto v(x,y;z)+y$ is nondecreasing for each $x\in \mathcal{I}_1$, $z\in \mathbb{R}^+$ (cf.~proof of Proposition \ref{prop:monot}) implies the second one.

\emph{3.} The proof of these two properties follows from the fact that $v(\cdot)$ is continuous by Proposition \ref{vcont} and Remark \ref{rem:cont}, and from point 2 above by using arguments as those employed in \cite[Prop.\ 2.2]{Jacka}.

\emph{4.} Notice that by \eqref{equivalences} one has
\beq
\label{equivalences2}
\{(x,y) \in \mathcal{I}_1 \times \mathcal{I}_2: z > z^*(x,y)\} = \{(x,y) \in \mathcal{I}_1 \times \mathcal{I}_2: v(x,y;z) > - y\},
\enq
for any $z \in \R^+$. The set on the right-hand side above is open since it is the preimage of an open set via the continuous mapping $(x,y) \mapsto v(x,y;z) + y$ (cf.\ Proposition \ref{vcont}). Hence the set on the left-hand side of \eqref{equivalences2} is open as well and thus $(x,y) \mapsto z^*(x,y)$ is upper-semicontinuous.\ep\\
\vspace{+3pt}

\noindent Now Proposition \ref{prop:z*} and the following
\begin{Assumption}\label{ass:cz}
$\displaystyle \lim_{z\uparrow\infty} c_z(x,z)=\infty$ for every $x\in\mathcal{I}_1$
\end{Assumption}
\noindent imply
\begin{Proposition}
\label{prop:z**}
Under Assumption \ref{ass:cz}, $\overline{z}$ is finite on $Q$.
\end{Proposition}
\noindent Then, thanks to Proposition \ref{prop:z*}-(1) one also has
\begin{Corollary}
\label{z*finite}
${z}^*$ is finite on $Q$.
\end{Corollary}

The topological characterization of the regions $\mathcal{C}$ and $\mathcal{A}$ is given in the following
\begin{Proposition}
$\mathcal{C}$ is open and $\mathcal{A}$ is closed. Moreover, under Assumption \ref{ass:cz}, they are connected.
\end{Proposition}
\textbf{Proof.}
The fact that $\mathcal{C}$ is open and $\mathcal{A}$ is closed follows from \eqref{regionsCA0} and Remark \ref{rem:cont}. Corollary \ref{z*finite} and \eqref{regionsCA1} imply the second part of the claim.\ep \\

\subsection{Optimal Control: a Verification Theorem}
The results obtained in Section \ref{sec:associated} on the optimal stopping problem \eqref{P:OS} 
allow us to provide the expression of the optimal control $\nu^*$ of problem \eqref{OCP} in terms of the  boundary $z^*$ of \eqref{defzstar}. Moreover, as a byproduct, we will also show that $V_z=v$ on $\Oc$ as expected (see Corollary \ref{cor:Vzv} below). 

Recall \eqref{P:OS} and define the functions
\beq
\label{Phi}
\Phi(x,z):=\mathbb{E}\bigg[\int_0^\infty e^{-rt} c(X_t^x,z)dt\bigg], \ \ \  (x,z)\in\mathcal{I}_1\times \mathbb{R}^+,
\enq
\beq
\label{varPhi}
\varphi(x,z):= \frac{\partial}{\partial z}\Phi(x,z)= \mathbb{E}\bigg[\int_0^\infty e^{-rt} c_z(X_t^x,z)dt\bigg], \ \ (x,z)\in\mathcal{I}_1\times \mathbb{R}^+,
\enq
and
\beq
\label{U}
U(x,y,z):= \Phi(x,z)-\int_z^\infty(v(x,y;q)-\varphi(x,q))dq, \ \ (x,y,z)\in \Oc.
\enq
Notice that $v(x,y;z)\geq\varphi(x,z)$ for every $(x,y,z)\in\Oc$, and therefore function $U$ in \eqref{U} above is well-defined (but, a priori, it may be equal to $-\infty$).

Introduce the nondecreasing process
\beq
\label{eq:optcont}
\nu^*_t:=\sup_{0\leq s\leq t}[z^*(X^x_s,Y^y_s)-z]^+, \quad t \geq 0, \quad \qquad \nu^*_{0-}=0,
\enq
with $z^*(x,y)$ as in \eqref{defzstar}. Notice that $\nu^*_t$ is the minimal amount of control needed at time $t \geq 0$ to keep $Z^{z,\nu^*}_t$ above $z^*(X^x_t,Y^y_t)$, thus solving a Skorokhod reflection problem.
\begin{Proposition}
\label{admissibilitynustar}
Under Assumption \ref{ass:cz} the process $\nu^*$ of \eqref{eq:optcont} is an admissible control.
\end{Proposition}
\textbf{Proof.}
Recall the set of admissible controls $\mathcal{V}$ of \eqref{setadmissblecontrols}. Clearly $\nu^*$ is a.s.\ finite thanks to Corollary \ref{z*finite}.
To prove that $\nu^* \in \mathcal{V}$ it remains to show that: \emph{i)} $t \mapsto \nu^*_t$ is right-continuous with left-limits; \emph{ii)} $\nu^*$ is ($\mathcal{F}_t$)-adapted.

We start by proving \emph{i)}. Clearly, $t\mapsto\nu^*_t$ admits left-limit at any point since it is nondecreasing. To show that $\nu^*$ has right-continuous sample paths, first notice that
\begin{align}
\label{usrc}
\limsup_{s \downarrow t}z^*(X^x_{s},Y^y_{s}) \leq z^*(X^x_t,Y^y_t)
\end{align}
by upper-semicontinuity of $z^*$ (cf.~Proposition \ref{prop:z*}) and continuity of $(X^x_{\cdot}, Y^y_{\cdot})$.
 Moreover, from \eqref{eq:optcont} and \eqref{usrc} we obtain
\begin{align}
\label{rc2}
  \lim_{s \downarrow t} \nu^*_s  =& \, \nu^*_t \vee \lim_{s \downarrow t} \sup_{t < u \leq s}[z^*(X^x_u,Y^y_u)-z]^+ \nonumber\\
  = & \, \nu^*_t \vee \limsup_{s \downarrow t} [z^*(X^x_s,Y^y_s)-z]^+ \leq \nu^*_t \vee [z^*(X^x_t,Y^y_t)-z]^+  =\nu^*_t.
\end{align}
Since $\lim_{s\downarrow t} \nu^*_s\ge \nu^*_t$ by monotonicity of $t\mapsto \nu^*_t$, then \eqref{rc2} implies right continuity.

As for \emph{ii)} the process $z^*(X^x,Y^y)$ is progressively measurable since it is the composition of the  Borel-measurable function $z^*$ (which is upper semicontinuous by Proposition \ref{prop:z*}) with the progressively measurable process $(X^x,Y^y)$. Therefore $\nu^*$ is progressively measurable by \cite[Th.\ IV.33, part (a)]{DellMeyerA}, hence adapted and \emph{ii)} above holds.
\ep

\begin{Theorem}
\label{teo:ver}
Let Assumption \ref{ass:cz} hold. Fix $(x,y,z)\in\Oc$  and take $\Phi(x,z)$, $\varphi(x,z)$ and $U(x,z)$ as in \eqref{Phi}, \eqref{varPhi} and \eqref{U}, respectively. Then one has $U(x,y,z) = V(x,y,z)$ and $\nu^*$ as in \eqref{eq:optcont} is optimal for the singular control problem \eqref{OCP}.
\end{Theorem}
\noindent It clearly follows from Theorem \ref{teo:ver} the following
\begin{Corollary}
\label{cor:Vzv}
The identity $V_z=v$ holds true on $\Oc$.
\end{Corollary}

The proof of Theorem \ref{teo:ver} is inspired by the arguments developed in \cite{BK} and \cite{KaratzasElKarouiSkorokhod}.
\medskip

\noindent\textbf{Proof of Theorem \ref{teo:ver}.} For $\nu \in \mathcal{V}$ define its right-continuous inverse (cf.\ \cite[Ch.\ 0, Sec.\ 4]{RY})
\beq
\label{taunu}
\tau^\nu(\xi):=\inf\{t\geq 0 \ | \ \nu_t > \xi\}, \qquad \xi\geq 0.
\enq
The process $\tau^\nu:=\{\tau^\nu(\xi),\  \xi \geq 0\}$ has increasing, right-continuous sample paths and hence it admits left-limits
\beq
\label{taunumeno}
\tau^\nu_{-}(\xi):=\inf\{t\geq 0 \ | \ \nu_t \geq \xi\}, \qquad \xi\geq 0.
\enq
The set of points $\xi\in\R^+$ at which $\tau^\nu(\xi)(\omega) \neq \tau^\nu_{-}(\xi)(\omega)$ is a.s.\ countable for a.e.~$\omega\in\Omega$.

Since $\nu$ is right-continuous and $\tau^\nu(\xi)$ is the first entry time of an open set, it is an $(\Fc_{t+})$-stopping time for any given and fixed $\xi \geq 0$. However, $(\Fc_t)_{t\ge0}$ is right-continuous (cf.~Section \ref{sec:IrrInvProb}), hence $\tau^\nu(\xi)$ is an $(\Fc_{t})$-stopping time. Moreover, $\tau^\nu_{-}(\xi)$ is the first entry time of the right-continuous process $\nu$ into a closed set and hence it is an $(\Fc_{t})$-stopping time as well for any $\xi \geq 0$. It then follows by the superharmonic characterization of $v$ that 
\begin{align}\label{ooi}
v(x,y;q)\geq \mathbb{E}\bigg[e^{-r \tau^{\nu}(\xi)}v(X^x_{\tau^{\nu}(\xi)}, Y^y_{\tau^{\nu}(\xi)};q)+\int_0^{\tau^{\nu}(\xi)}e^{-rs}c_z(X_s^x,q)ds\bigg],
\end{align}
for any $\xi \geq 0$ and $(x,y,q)\in\Oc$. Then, for any $(x,y,z)\in\Oc$, taking $\xi=q-z$, $q\ge z$ in \eqref{ooi} and recalling \eqref{Phi}, \eqref{varPhi}, and \eqref{U}, we obtain
\begin{eqnarray}
\label{verifico1}
U(x,y,z)-\Phi(x,z)& \hspace{-0.25cm}\leq &\hspace{-0.25cm} -\int_z^\infty\bigg(\mathbb{E}\bigg[e^{-r\tau^{\nu}(q-z)}v(X^x_{\tau^{\nu}(q-z)}, Y^y_{\tau^{\nu}(q-z)};q)+ \nonumber \\
&& + \int_0^{\tau^{\nu}(q-z)}e^{-rs}c_z(X_s^x,q)ds\bigg]\bigg)dq + \int_z^\infty\mathbb{E}\bigg[\int_0^\infty e^{-rs}c_z(X^x_s,q)ds\bigg]dq  \nonumber \\
&\hspace{-0.25cm} \leq \hspace{-0.25cm} &\int_z^\infty \mathbb{E}\bigg[e^{-r\tau^{\nu}(q-z)}Y_{\tau^{\nu}(q-z)}^y\bigg]dq
-\int_z^\infty\mathbb{E}\bigg[\int_0^{\tau^{\nu}(q-z)}e^{-rs}c_z(X_s^x,q)ds\bigg]dq \nonumber \\
&& \hspace{2cm} +\int_z^\infty\mathbb{E}\bigg[\int_0^\infty e^{-rs}c_z(X^x_s,q)ds\bigg]dq,
\end{eqnarray}
where we have used that $v(\,\cdot\,,\zeta;\,\cdot\,)\geq -\zeta$ (cf.\ Proposition \ref{prop:v}) in the second inequality.
We now claim (and we will prove it later) that we can apply Fubini-Tonelli's  Theorem in the last expression of \eqref{verifico1} to obtain
\begin{eqnarray}
\label{verifico2}
U(x,y,z)-\Phi(x,z)& \hspace{-0.25cm} \leq &\hspace{-0.25cm} \mathbb{E}\bigg[\int_z^\infty e^{-r\tau^{\nu}(q-z)}Y_{\tau^{\nu}(q-z)}^y dq
-\int_z^\infty\bigg(\int_0^{\tau^{\nu}(q-z)} e^{-rs}c_z(X_s^x,q)ds\bigg)\,dq\bigg] \nonumber \\
& & \hspace{2cm} +\, \mathbb{E}\bigg[\int_z^\infty \bigg(\int_0^\infty e^{-rs} c_z(X^x_s,q)ds\bigg)\,dq\bigg].
\end{eqnarray}
The change of variable formula of \cite[Ch.\ 0, Prop.\ 4.9]{RY} (see also \cite[eq.~(4.7)]{BK}) implies
\beq
\label{changevariable}
\int_z^\infty e^{-r\tau^{\nu}(q-z)}Y_{\tau^{\nu}(q-z)}^y dq = \int_0^\infty e^{-rs} Y^y_s d\nu_s.
\enq
Moreover, $\tau^{\nu}(q-z) < s$ if and only if $\nu_s > q - z$, where $s \geq 0$. Therefore, from \eqref{verifico2} and \eqref{changevariable} we obtain
\begin{eqnarray}
\label{verifico3}
U(x,y,z)-\Phi(x,z) & \hspace{-0.25cm} \leq \hspace{-0.25cm} & \mathbb{E}\bigg[ \int_0^\infty e^{-rs} Y^y_sd\nu_s
+\int_z^\infty \bigg(\int_{\tau^{\nu}(q-z)}^\infty e^{-rs}c_z(X_s^x,q)ds\bigg)dq\bigg] \nonumber \\
& \hspace{-0.25cm} = \hspace{-0.25cm} & \mathbb{E}\bigg[ \int_0^\infty e^{-rs} Y^y_sd\nu_s
+\int_z^\infty \bigg(\int_{0}^\infty e^{-rs}c_z(X_s^x,q)\mathbf{1}_{\{\nu_s > q - z\}}ds\bigg)dq\bigg] \nonumber \\
& \hspace{-0.25cm} = \hspace{-0.25cm} & \mathbb{E}\bigg[ \int_0^\infty e^{-rs} Y^y_sd\nu_s
+\int_{0}^\infty e^{-rs}\bigg(\int_z^{z+\nu_s}  c_z(X_s^x,q)dq\bigg)ds\bigg] \\
& \hspace{-0.25cm} = \hspace{-0.25cm} & \mathbb{E}\bigg[ \int_0^\infty e^{-rs} Y^y_sd\nu_s
+\int_{0}^\infty e^{-rs} \Big[c(X_s^x,Z_s^{z,\nu})-c(X^x_s,z)\Big]ds\bigg] \nonumber \\
& \hspace{-0.25cm} = \hspace{-0.25cm} & \mathcal{J}_{x,y,z}(\nu)-\Phi(x,z). \nonumber
\end{eqnarray}
Since $\nu\in\mathcal{V}$ is arbitrary, it follows
\beq\label{ineq}
U(x,y,z)\leq V(x,y,z).
\enq
Now we want to show that picking $\nu^*$ as in \eqref{eq:optcont} in the arguments above all the inequalities become equalities due to \eqref{superarmonic2}. First, notice that \eqref{superarmonic2}, \eqref{regionsCA0}, and \eqref{regionsCA1} give
\beq\label{equiv1}
\tau^*(x,y;q)= \ \inf\{t\geq 0 \ | \ z^*(X^x_t,Y^y_t) \geq q \}.
\enq
Then, fix $z\in\R^+$, take $t\ge 0$ arbitrary, and note that, by \eqref{taunumeno} and \eqref{equiv1}, we have $\P$-a.s.~the equivalences
\begin{align*}
 &\tau^{\nu^*}_{-}(q-z) \leq t \ \Longleftrightarrow\  \nu^*_t \geq q - z\  \Longleftrightarrow\   \sup_{0\leq s\leq t}[z^*(X^x_s,Y^y_s)-z]^+ \geq q-z \\
 &\Longleftrightarrow \  z^*(X^x_{\theta},Y^y_{\theta}) \geq q \ \mbox{ for some } \theta \in [0,t] \  \Longleftrightarrow \ \tau^*(x,y;q) \leq t.
\end{align*}
So, we can conclude that $\tau^{\nu^*}_-(q-z)=\tau^*(x,y;q)$ $\P$-a.s.~and for a.e.~$q\ge z$. However, by \eqref{taunu} and \eqref{taunumeno}, we also have $\tau^{\nu^*}_-(q-z)=\tau^{\nu^*}(q-z)$ $\P$-a.s.~and for a.e.~$q\ge z$; hence
\begin{align}\label{equivtda}
\text{$\tau^{\nu^*}(q-z)=\tau^*(x,y;q)$ $\P$-a.s.~and for a.e.~$q\ge z$.}
\end{align}

Now take $\nu=\nu^*$ and $\xi=q-z$ in order to obtain equality in \eqref{ooi}, by harmonic property of $v$ in the continuation set. 
Optimality of $\tau^*=\tau^{\nu^*}$ (cf.~\eqref{equivtda}) also gives equality in \eqref{verifico1}; then, we can interchange the integrals and argue as in \eqref{verifico2} and \eqref{verifico3} to obtain $U(x,y,z)=\mathcal{J}_{x,y,z}(\nu^*)$. Then $U=V$ on $\Oc$, by \eqref{ineq}, and $\nu^*$ is optimal.

To conclude the proof, we need to show that we could actually interchange the order of integration in \eqref{verifico1} to get \eqref{verifico2}.
Clearly
$$\int_z^\infty \mathbb{E}\bigg[e^{-r\tau^{\nu}(q-z)}Y_{\tau^{\nu}(q-z)}^y\bigg]dq = \mathbb{E}\bigg[\int_z^\infty e^{-r\tau^{\nu}(q-z)}Y_{\tau^{\nu}(q-z)}^y dq\bigg],$$
by Tonelli's Theorem, since $Y^y$ has positive sample paths. Therefore, we have only to show that
\beq
\label{checkFubini}
\mathbb{E}\bigg[\int_z^\infty \Big(\int_{\tau^{\nu}(q-z)}^{\infty}e^{-rs}|c_z(X_s^x,q)|ds\Big) dq\bigg] < \infty.
\enq
Define
$$q^*_s:=\inf\{q \in \R: c_z(X^x_s,q) > 0\},$$
which exists and is unique, since $c(x,\cdot)$ is convex. Now, recall that $\tau^{\nu}(q-z) < s$ if and only if $\nu_s > q - z$, $s \geq 0$, and notice that any admissible control $\nu$ should also satisfy, without loss of generality,
\begin{align}\label{adm:int}
\E\bigg[\int_0^\infty{e^{-r\,t}c(X^x_t,z+\nu_t)dt}\bigg]<+\infty.
\end{align}
Indeed, \eqref{adm:int} holds for the optimal control $\nu^*$ (if it exists), since $J_{x,y,z}(\nu^*)\le J_{x,y,z}(0)$.
Then, Tonelli's Theorem, \eqref{adm:int}, and the fact that $c\ge 0$ give
\begin{eqnarray*}
& & \mathbb{E}\bigg[\int_z^\infty \Big(\int_{\tau^{\nu}(q-z)}^{\infty}e^{-rs}|c_z(X_s^x,q)|ds\Big) dq\bigg] = \mathbb{E}\bigg[\int_z^\infty \Big(\int_0^{\infty}e^{-rs}|c_z(X_s^x,q)|\mathds{1}_{\{\tau^{\nu}(q-z) < s\}}ds\Big) dq\bigg] \nonumber \\
& & = \mathbb{E}\bigg[\int_0^{\infty}e^{-rs} \Big(\int_z^{z+\nu_s}|c_z(X_s^x,q)|dq\Big) ds\bigg] = \mathbb{E}\bigg[\int_0^{\infty}e^{-rs} \Big(\int_{(z+\nu_s) \wedge q^*_s}^{z+\nu_s}c_z(X_s^x,q)dq\Big) ds\bigg] \nonumber \\
& & \hspace{1cm} - \mathbb{E}\bigg[\int_0^{\infty}e^{-rs} \Big(\int_{z}^{(z+\nu_s) \wedge q^*_s}c_z(X_s^x,q)dq\Big) ds\bigg] \nonumber \\
& & \leq  \mathbb{E}\bigg[\int_0^{\infty}e^{-rs} c(X_s^x,z) ds + \int_0^{\infty}e^{-rs} c(X_s^x,z+\nu_s) ds\bigg] < \infty. \nonumber
\end{eqnarray*}
\ep


\appendix

\section{Appendix}
\label{app:os1d}
\renewcommand{\theequation}{A-\arabic{equation}}

\subsection{Proof of Proposition \ref{teo:ex}}
\label{ProofProp313}

\emph{Step 1.}\hspace{0.15cm} Since $\mu_i$, $\sigma_i$, $i=1,2$ are bounded and continuous on $Q_n$, existence and uniqueness of a function $u_n(\cdot\,;z)\in W^{2,p}(Q_n)$ for all $1\leq p<\infty$ solving \eqref{VIn} in the a.e.\ sense and satisfying \eqref{BC} follow by \cite[Ch.\,I, Th.\,3.2 and Th.\,3.4]{Fri}. The function $u_n(\,\cdot\,;z)$ can be continuously extended outside $Q_n$ by setting
\beq
\label{extension}
u_n(x,y;z)=-y,\quad \quad (x,y)\in Q\setminus Q_n,
\enq
and we denote such extension again by $u_n$, with a slight abuse of notation.
\smallskip

\emph{Step 2.} We now show that $v_n(\,\cdot\,;z) = u_n(\,\cdot\,;z)$ over $Q_n$ and that the stopping time \eqref{opt-st-n} is optimal for problem \eqref{P:OSn}.

If $(x,y)\in Q\setminus Q_n$, then the claim clearly follows from Proposition \ref{prop:vR}-(2).
Assume $(x,y)\in Q_n$; since $u_n\in W^{2,p}(Q_n)$, by \cite[Ch.~7.6]{GT} we can find a sequence $\big\{u^{\,k}_n(\,\cdot\,;z)\,,k\in\mathbb{N}\,\big\}\subset C^\infty(Q)$ such that ${u}^{\,k}_n (\,\cdot\,;z)\rightarrow {u}_n(\,\cdot\,;z)$ in $W^{2,p}(Q_n)$, $p\in[1,+\infty)$, as $k\rightarrow \infty$. Moreover, since ${u}_n$ is continuous and $\overline{Q}_n$ is a compact, we have  ${u}^{\,k}_n(\,\cdot\,;z)\rightarrow {u}_n(\,\cdot\,;z)$ uniformly on $\overline{Q}_n$ (cf.\ \cite[Ch.~7.2, Lemma 7.1]{GT}).

Dynkin's formula yields for any bounded stopping time $\tau$
\begin{align}\label{dyn041}
u^k_n(x,y;z)=\E\left[e^{-r(\tau\wedge\sigma_n)}u_n^k(X^x_{\tau\wedge\sigma_n},
Y^y_{\tau\wedge\sigma_n};z)-\int^{\tau\wedge\sigma_n}_0{e^{-rt}(\mathbb{L}-r)u_n^k(X^x_t,Y^y_t;z)\,dt}
\right].
\end{align}
Then, by localization arguments and using \eqref{conventions1}, we conclude that \eqref{dyn041} actually holds for any $\tau \in \mathcal{T}$.
We claim (and we will prove it later) that taking limits as $k\rightarrow \infty$ in \eqref{dyn041} leads to
{\small{\begin{align}
\label{dyn02}
{u}_n(x,y;z)=\E\left[e^{-r(\tau\wedge\sigma_n)}{u}_n(X^x_{\tau\wedge\sigma_n},
Y^y_{\tau\wedge\sigma_n};z)-\int^{\tau\wedge\sigma_n}_0{e^{-rt}(\mathbb{L}-r){u}_n(X^x_t,Y^y_t;z)\,dt}
\right], \  \forall\tau \in \mathcal{T}.
\end{align}}}
The right-hand side of \eqref{dyn02} is well defined, since Assumption \ref{ass:laws} implies that the law of $(X^x,Y^y)$ is absolutely continuous with respect to the Lebesgue measure and $(\mathbb{L}-r)u_n$ is defined up to a Lebesgue null-measure set. We now use the variational inequality \eqref{VIn} in \eqref{dyn02} to obtain
\begin{align}\label{dyn05}
u_n(x,y;z)\ge\E\left[-e^{-r(\tau\wedge\sigma_n)}Y^y_{\tau\wedge\sigma_n}+\int^{\tau\wedge\sigma_n}_0 e^{-rt}c_z(X^x_t,z)\,dt
\right].
\end{align}
Hence, by arbitrariness of $\tau$, one has $u_n(x,y;z)\ge v_n(x,y;z)$.

To obtain the reverse inequality, take
\beq
\label{taucandidate}
\tau:=\inf\big\{t\ge0\,|\,u_n(X^x_t,Y^y_t;z)=-Y^y_t\big\}
\enq
in \eqref{dyn02} and recall that $u_n=-y$ on $Q\setminus Q_n$, that $u_n\in C^0(\overline{Q}_n)$ (cf.~Remark \ref{Sobolev}), and that  $\overline{Q}_n$ is bounded, so that $u_n$ is bounded in $\overline{Q}_n$ as well. It follows that
\begin{align}\label{zero}
e^{-r(\tau\wedge\sigma_n)}u_n(X^x_{\tau\wedge\sigma_n},Y^y_{\tau\wedge\sigma_n};z)=&
\,e^{-r(\tau\wedge\sigma_n)}u_n(X^x_{\tau\wedge\sigma_n},Y^y_{\tau\wedge\sigma_n};z)\mathds{1}_{\{\tau\wedge\sigma_n<\infty\}}\nonumber\\
=&-e^{-r(\tau\wedge\sigma_n)}Y^y_{\tau\wedge\sigma_n}\mathds{1}_{\{\tau\wedge\sigma_n<\infty\}}=
-e^{-r(\tau\wedge\sigma_n)}Y^y_{\tau\wedge\sigma_n}\quad\text{$\P$-a.s.}
\end{align}
by \eqref{conventions} and \eqref{conventions1}. Moreover,  by \eqref{VIn}, we have $(\mathbb{L}_X-r)u_n=-c_z$ on the set $\big\{(x,y)\in Q_n\,|\,u_n(x,y;z)>-y\big\}$. Hence, \eqref{dyn02} and \eqref{zero} give
\begin{align}
u_n(x,y;z)=\E\left[-e^{-r(\tau\wedge\sigma_n)}Y^y_{\tau\wedge\sigma_n}+\int^{\tau\wedge\sigma_n}_0 e^{-rt}c_z(X^x_t,z)\,dt
\right]\le v_n(x,y;z).
\end{align}
Therefore, we conclude that $u_n=v_n$ on $Q$, and that the stopping time $\tau$ defined in \eqref{taucandidate} is optimal for problem \eqref{P:OSn} and coincides with the stopping time $\tau^*_n(x,y;z)$ defined in \eqref{opt-st-n}.

Now, to complete the proof, we only need to show that \eqref{dyn02} follows from \eqref{dyn041} as $k \rightarrow  \infty$.
In fact, the term on the left-hand side of \eqref{dyn041} converges pointwisely and the first term in the expectation on the right-hand side converges by uniform convergence. To check convergence of the integral term in the expectation on the right-hand side, we take $q_n>1$ as in Assumption \ref{ass:laws}-(2), $p_n$ such that $\frac{1}{p_n}+\frac{1}{q_n}=1$, and, for simplicity, denote $q:=q_n$ and $p:=p_n$. Then, by H\"older's inequality, we have
\begin{align}\label{dyn03}
&\left|\,\E\left[\int^{\tau\wedge\sigma_n}_0\hspace{-8pt}{e^{-rt}(\mathbb{L}-r)({u}^{\,k}_n-{u}_n)
(X^x_t,Y^y_t;z)\,dt}
\right]\,\right|\nonumber \\
&\le\int^\infty_0{e^{-rt}\left(\int_{Q_n}\big|(\mathbb{L}-r)
({u}^{\,k}_n-{u}_n)(\xi,\zeta;z)\big|\,
p_1(t,x,\xi)p_2(t,y,\zeta)d\xi\,d\zeta\right)dt}\\ &\le C_{M_1,M_2,r,n}\,\big\|{u}^{\,k}_n-{u}_n\big\|_{W^{2,p}(Q_n)},\nonumber
\end{align}
where last inequality  follows by Assumptions \ref{ass:D2}-(i) and \ref{ass:laws}-(2), with $C_{M_1,M_2,r,n}>0$ depending on $Q_n$,\,$r$, and  $M_i:=\sup_{\overline{Q}_n}\left\{ |\mu_i|+|\sigma_i|\right\}$\,, $i=1,2$. Now, the right-hand side of \eqref{dyn03} vanishes as $k\rightarrow \infty$ by definition of $u^k_n$.
\ep

\subsection{Two Technical Lemmas}
\label{LemmaAppendix}

$1.$ By \eqref{def:CnAn}, heuristically one has $(\mathbb{L}-r)v_n(x,y)=ry-\mu_2(y)$ on $\mathcal{A}^n_z\cap Q_n$. However, at this stage we do not have sufficient information about the topological properties of $\mathcal{A}^n_z$ (for example it could have positive measure, but a priori also empty interior part). The following Lemma provides a rigorous statement and proof of the previous equality. 
\begin{Lemma}
\label{Lemmahessiano}
One has
\begin{equation}\label{operatore}
(\mathbb{L}-r)v_n(x,y)=ry-\mu_2(y), \ \ \ \ \mbox{for a.e.} \ (x,y)\in \mathcal{A}^n_z\cap Q_n.
\end{equation}
\end{Lemma}
\textbf{Proof.}
Recall that $v_n(\,\cdot\,;z) \in W^{2,p}(Q_n)$ for any $p \in [1,\infty)$ (cf.\ Proposition \ref{teo:ex}). Set $\bar{v}_n(x,y;z):=v_n(x,y;z)+y$, hence $\bar{v}_n\in C^1(Q_n)$ by Sobolev's embedding (see for instance \cite[Ch.\ 9, Cor.\ 9.15]{Br}) and proving \eqref{operatore} amounts to showing that $(\mathbb{L}-r) \bar{v}_n = 0$ a.e.\ on $\mathcal{A}^n_z\cap Q_n$. Since $\bar{v}_n=0$ over $\mathcal{A}^n_z$, it must also be $\nabla\bar{v}_n=0$ over $\mathcal{A}^n_z\cap Q_n$.
To complete the proof, it thus remains to show that the Hessian matrix $\mathcal{A}^n_z$; that is,
$D^2\bar{v}_n$ is zero a.e.\ over $\mathcal{A}^n_z\cap Q_n$.
This follows by \cite[Cor.\ 1-(i), p.\ 84]{EG}\footnote{It is worth noting that \cite[Cor.\ 1-(i), p.\ 84]{EG} requires $f$ to be Lipschitz continuous, which is not guaranteed for us. However, Lipschitz continuity is only needed there to have existence a.e.\ of the gradient $\nabla f$,
which we have due to \cite[Th.\ 1, p.\ 235]{EG}, since $\nabla \bar{v}_n\in W^{1,p}(Q_n)$.}, with $f$ therein defined by $f:=\nabla \bar{v}_n$.
\ep
\vspace{+8pt}

$2.$ The next result is important for the proof of Theorem \ref{teo:uniqueness}.
\begin{Lemma}\label{lemma:w}
Let Assumptions \ref{ass:mu21}, \ref{ass:densitiesderivative}, \ref{ass:mu2} hold, and assume that  $\mathcal{C}_z\neq \emptyset$ and $\mathcal{A}_z\neq \emptyset$. Let $\hat{y}(\,\cdot\,;z):\mathcal{I}_1\rightarrow\overline{\mathcal{I}}_2$ be a solution of \eqref{eq:intan2} and take $w$ as in \eqref{defW}. Then $v(\,\cdot\,;z)\geq w(\,\cdot\,;z)$ on $Q$.
\end{Lemma}
\textbf{Proof.} Recall the notation introduced in \eqref{notationhat}.

\emph{Step 1.} Since $\hat{y}(\,\cdot\,;z)$ is a solution of \eqref{eq:intan2}, i.e.\ of \eqref{eq:int}, it is easy to see that
 $w$ of \eqref{defW} verifies
\begin{align}\label{eq:ww00}
w(x,\hat{y}(x;z);z)=-\hat{y}(x;z),\qquad\forall x\in\hat{\mathcal{D}}_z,
\end{align}
 therefore
\begin{align}\label{eq:ww001}
w(x,\hat{y}(x;z);z)\leq v(x,\hat{y}(x;z);z),\qquad\forall x\in\hat{\mathcal{D}}_z,
\end{align}

\smallskip
\emph{Step 2.}
Here we show that
\begin{align}\label{ww}
w(x,y;z)=-y, \ \ \ \ \forall y<\hat{y}(x;z),  \ \forall x\in\hat{\Dc}_z\cup [\hat{x},\overline{x}),
\end{align}
which implies
$$w(x,y;z)\le v(x,y;z),\ \ \ \forall y<\hat{y}(x;z),   \ x\in\hat{\Dc}_z\cup [\hat{x},\overline{x}).$$
Take $x\in\hat{\Dc}_z\cup [\hat{x},\overline{x})$, $y<\hat{y}(x;z)$ and define $\sigma = {{\sigma}}(x,y;z):= \inf\big\{t\ge0\,|\,Y^y_t\ge \hat{y}(X^x_t;z)\big\}$. By definition of $\hat{y}(\,\cdot\,;z)$ and ${\sigma}$, we have
\begin{align}\label{www01}
\widehat{H}(X^x_t,Y^y_t;z)=-\big(rY^y_t-\mu_2(Y^y_t)\big), \ \ \text{$\forall t\le {\sigma}$, $\P$-a.s.}
\end{align}
Then, using the martingale property \eqref{eq:mg02} up to the stopping time $\sigma\wedge n$, $n\in\mathbb{N}$, it follows by \eqref{www01} that
\begin{align}\label{www00}
w(x,y;z)=\E\bigg[-e^{-r\sigma}Y^y_{{\sigma}}\mathds{1}_{\{\sigma\le n\}}+e^{-r n}w(X^x_{n},Y^y_{n};z)\mathds{1}_{\{\sigma>n\}}-\int^{\sigma\wedge n}_0{e^{-rt}\big(rY^y_t-\mu_2(Y^y_t)\big)dt}\bigg].
\end{align}
Assumption \ref{Ass:r}, \eqref{uniq02}, and the bound \eqref{eq:subpon02} give in the limit as $n\rightarrow\infty$
\begin{align}\label{www02}
w(x,y;z)=\E\bigg[-e^{-r\sigma}Y^y_{{\sigma}}-\int^{\sigma}_0{e^{-rt}\big(rY^y_t-\mu_2(Y^y_t)\big)dt}\bigg]=-y,
\end{align}
where the last equality follows by Lemma \ref{ibp}. Hence, \eqref{ww} is proved.

\smallskip
\emph{Step 3.}
Here we prove that
\beq\label{vw1}
w(x,y;z)\leq v(x,y;z), \ \ \ \ \forall y>\hat{y}(x;z),  \ \forall x\in(\underline{x},\check{x}]\cup\hat{\Dc}_z.
\enq
Take $x\in (\underline{x},\check{x}]\cup\hat{\Dc}_z$, $y>\hat{y}(x;z)$, and consider the stopping time
$${\tau}={\tau}(x,y;z):=\inf\big\{t\ge0\ |\ Y^y_t\le \hat{y}(X^x_t;z)\big\}.$$
By definitions of $\hat{y}(\,\cdot\,;z)$ and ${\tau}$, and by using the same localization argument as in \emph{Step 2} above, we obtain
\begin{align}\label{uniq03}
w(x,y;z)=\E\left[-e^{-r{{\tau}}}Y^y_{{\tau}}+\int_0^{{\tau}}{e^{-rs}c_z(X^x_t,z)dt}\right]\le v(x,y;z).
\end{align}

\emph{Step 4.} Now Lemma \ref{lemma:w} follows by \eqref{eq:ww00}, \eqref{ww}, and \eqref{vw1}.
\ep

\subsection{Some properties of non-entrance boundaries}
\label{A1}

Here we establish some properties of the diffusion $Y^y$ having natural boundaries (cf.\ Assumption \ref{ass:D2}), hence \emph{non-entrance}. 
We prove the following results for the lower boundary $\underline{y}$, but similar arguments also hold for $\overline{y}$ if it is finite.

By definition of a non-entrance boundary (see, e.g., \cite[p.\,305]{RY}) we have 
\begin{equation}
\label{ppp}
\lim_{y\downarrow\underline{y}}\mathbb{P}_{y}\{\tau_z<t\}=0, \ \ \ \forall z>\underline y, \ t>0,
\end{equation}
where
$\tau_z:=\inf\{ s\geq 0 \ | \ Y^y_s\geq z\}$.
Taking an arbitrary $\varepsilon>0$, a given and fixed $t>0$, and setting $z:=z_{\varepsilon}=\underline y+\varepsilon$, we have
\begin{equation}
\label{ppp-bis}
\{|Y_t^y-\underline y| > \varepsilon\}\subseteq \{\sup_{s\in[0,t]} Y^y_s> z\}= \{\tau_{z}<t\}.
\end{equation}
It thus follows from \eqref{ppp-bis} and \eqref{ppp} that
$Y_t^y\rightarrow \underline y$ in probability (hence in law) as $y\downarrow \underline y$ for every $t>0$ given and fixed; that is, the  probability measure on $\Ic_2$ with density $p_2(t,y,\cdot)$, (cf.\ Assumption \ref{ass:laws}) converges weakly to the Dirac's delta measure $\delta_{\underline{y}}(\cdot)$,  when $y \downarrow \underline{y}$.
Therefore, by dominated convergence, one also has
\begin{align}\label{mandl}
\lim_{y\downarrow \underline{y}} \E\bigg[\int^\infty_0{e^{-rt}f(Y^y_t)dt}\bigg]=\frac{1}{r}f(\underline{y})\qquad\forall f\in C_b(\R).
\end{align}
\medskip

We now show that $\mu_2(\underline{y})=\sigma_2(\underline{y})=0$. The same holds for $\overline{y}$ if it is finite. 
\medskip

\emph{Case 1.} If $\Ic_2$ is bounded, an application of Dynkin's formula to any $g\in C^2_b(\R)$ leads to
\begin{align}
g(y)=-\E\bigg[\int^\infty_0{e^{-rt}\Big(\frac{1}{2}\sigma^2_2(Y^y_t)g''(Y^y_t)+\mu_2(Y^y_t)g'(Y^y_t)-rg(Y^y_t)\Big)dt}\bigg].
\end{align}
Then taking limits as $y\downarrow \underline{y}$, noting that $\mu_2$ and $\sigma_2$ are bounded and continuous and by applying \eqref{mandl} we get
\begin{align}
\frac{1}{2}\sigma^2_2(\underline{y})g''(\underline{y})+\mu_2(\underline{y})g'(\underline{y})=0,
\end{align}
and since $g$ is arbitrary it must be $\mu_2(\underline{y})=\sigma_2(\underline{y})=0$.
\medskip

\emph{Case 2}. If $\Ic_2$ is unbounded (i.e.\ if $\Ic_2=(\underline{y},\infty)$), we approximate $(\mu_2,\sigma_2)$ by continuous bounded functions $(\mu^n_2,\sigma^n_2)$ such that $\mu^n_2=\mu_2$ and $\sigma^n_2=\sigma_2$ on $[\underline{y},n\vee\underline{y}]$ with $\mu^n_2(y)\rightarrow\mu_2(y)$ and $\sigma_2^n(y)\rightarrow\sigma_2(y)$ as $n\rightarrow \infty$ pointwise on $\Ic_2$. For $y\in(\underline{y},n\vee\underline{y})$ the associated diffusion with coefficients $\mu^n_2$ and $\sigma^n_2$, denoted by $Y^{y,n}$, coincides with $Y^{y}$ up to the first exit time from $(\underline{y},n\vee\underline{y})$ by uniqueness of the solution of \eqref{state:Y}; moreover, $\underline{y}$ is a natural boundary for $Y^{y,n}$ as well.
Repeating arguments as in \emph{Case 1} above we get $\mu_2^n(\underline{y})=\sigma_2^n(\underline{y})=0$ for all $n\in\N$, thus $\mu_2(\underline{y})=\sigma_2(\underline{y})=0$.


%
\bigskip

\subsection{Discussion on Problem \eqref{eq:undv}}
\label{A2}

Problem \eqref{eq:undv} is standard in the optimal stopping literature (cf.~for instance \cite{PeskShir} for methods of solution) and hence we only sketch arguments leading to its main properties. It is easy to see that $x\mapsto \underline{v}(x;z)$ is nonincreasing and hence there exists $b_*\in\overline{\Ic}_1$ such that $\underline{\Ac}_z=[b_*,\overline{x})$, where the boundary value $\overline{x}$ cannot be included as otherwise $\Ac_z=\emptyset$ thus contradicting the assumption of Proposition \ref{prop:shapeD}. It is possible to show that $\underline{v}(\,\cdot\,;z)\in C^1(\Ic_1)$, $\underline{v}_{xx}(\,\cdot\,;z)$ is locally bounded at $b_*$ and hence that the probabilistic representation
\begin{align}
\label{rappresentationvsotto}
\underline{v}(x;z)=\E\Big[\int_0^\infty{e^{-rt}\Big(c_z(X^x_t;z)\mathds{1}_{\{X^x_t < b_*\}}-r\underline{y}\mathds{1}_{\{X^x_t\ge b_*\}}\Big)dt}\Big]
\end{align}
holds by It\^o-Tanaka formula.
Since \eqref{rappresentationvsotto} holds for any $x \in \Ic_1$, then if $b_*\in\Ic_1$ by evaluating \eqref{rappresentationvsotto} for $x = b_*$, one easily finds that $b_*$ solves \eqref{eq:x*1}. Arguments similar to (but simpler than) those employed in the proof of Theorem \ref{teo:uniqueness} show that \eqref{eq:x*1} admits a unique solution in $(\theta_*, \overline{x})$ and therefore it must be $\tilde{x}=b_*$.
On the other hand, if $b_* = \underline{x}$, repeating arguments as those of the proof of Theorem \ref{teo:uniqueness}, \emph{Step 2}, one can show that $\tilde{x}=b_*$, thus concluding.
\bigskip

\textbf{Acknowledgments.} The authors thank two anonymous referees for their pertinent comments, and Goran Peskir, Frank Riedel and Mauro Rosestolato for useful suggestions and references.


\end{document}